\documentclass[11pt]{article}
\usepackage{latexsym,amssymb}

\renewcommand{\baselinestretch}{1.20}
\newcommand{\C}[1]{{\protect\cal #1}}
\newcommand{\B}[1]{{\bf #1}}
\newcommand{\I}[1]{{\mathbb #1}}
\newcommand{\OO}[1]{\overline{#1}}
\newcommand{\e}{\varepsilon}
\renewcommand{\mid}{:}

\newcommand{\beq}[1]{\begin{equation}\label{eq:#1}}
\newcommand{\eeq}{\end{equation}}
\newcommand{\req}[1]{\textrm{(\ref{eq:#1})}}

\newtheorem{theorem}{Theorem}
\newtheorem{lemma}[theorem]{Lemma}

\newtheorem{proposition}[theorem]{Proposition}

\newtheorem{corollary}[theorem]{Corollary}
\newtheorem{question}[theorem]{Question}

\newtheorem{cla}{Claim}[theorem]
\newcommand{\bcl}[2]{\begin{cla}\label{cl:#1}#2\end{cla}}

\newcommand{\bpf}[1][Proof.]{\smallskip\noindent{\it #1} }
\newcommand{\qed}{\nolinebreak\mbox{\hspace{5 true pt}%
  \rule[-0.85 true pt]{3.9 true pt}{8.1 true pt}}}
\newcommand{\cqed}{\nolinebreak\mbox{\hspace{5 true pt}%
  \rule[-0.85 true pt]{2.0 true pt}{8.1 true pt}}}
\newcommand{\epf}{\qed \medskip}
\newcommand{\bcpf}{\bpf[Proof of Claim.]}
\newcommand{\ecpf}{\cqed \medskip}

\newcommand{\OPdata}{Oleg Pikhurko\footnote{The author was supported by the
European Research Council (grant agreement no.~306493)
and the National Science Foundation of the USA (grant DMS-1100215). This project
was initiated during the workshop \emph{``Hypergraph Tur\'an Problem''}
held at the American Institute of Mathematics, Palo Alto, 
March 21--25, 2011.
}\\
Mathematics Institute and DIMAP\\
University of Warwick\\
Coventry CV4 7AL, UK}
\date{}

\newcommand{\IS}{\I S} 
\newcommand{\ex}{\mathrm{ex}}
\newcommand{\PI}[2]{\Pi^{(#1)}_{#2}}
\newcommand{\FinPi}[1]{\PI{#1}{\mathrm{fin}}}
\newcommand{\InfPi}[1]{\Pi_{\infty}^{(#1)}}
\newcommand{\branch}[2]{\mathrm{br}_{#1}(#2)}
\newcommand{\blow}[2]{#1(\!(#2)\!)}
\newcommand{\density}[1]{\Lambda_{#1}}
\newcommand{\multiset}[1]{\{\hspace{-0.25em}\{\hspace{0.1em}#1\hspace{0.1em}\}\hspace{-0.25em}\}}

\newcommand{\LIM}[1]{\mathrm{LIM}^{(#1)}}
\newcommand{\CG}{\C G}
\newcommand{\cl}[1]{\OO{#1}}
\newcommand{\rep}[2]{{#1}^{(#2)}}
\newcommand{\hide}[1]{}

\begin{document}

\author{\OPdata}
\title{On Possible Tur\'an Densities}
\maketitle

\begin{abstract}
 The \emph{Tur\'an density} $\pi(\C F)$ of a family $\C F$ of
$k$-graphs is the limit as $n\to\infty$ 
of the maximum edge density 
of an $\C F$-free $k$-graph on $n$ vertices.
 Let $\InfPi{k}$ consist of all possible Tur\'an densities and
let $\FinPi{k}\subseteq \InfPi{k}$ be the set of Tur\'an densities of finite $k$-graph families.

Here we prove that $\FinPi{k}$ contains 
every density obtained from an arbitrary
finite construction by optimally blowing it up and using
recursion inside the specified set of parts. As an application, 
we show that $\FinPi{k}$ contains an irrational
number for each $k\ge 3$.

Also, we show that $\InfPi{k}$ has cardinality of the continuum.
In particular, $\InfPi{k}\not=\FinPi{k}$.
\end{abstract}

\section{Introduction}\label{Intro}

Let $\C F$ be a (possibly infinite) family of \emph{$k$-graphs} (that is,
$k$-uniform set systems). We call elements of $\C F$ \emph{forbidden}. A
$k$-graph $G$ is \emph{$\C F$-free} if no member $F\in \C F$ is 
a subgraph of $G$, that is, we cannot obtain $F$
by deleting some vertices and edges from $G$. The \emph{Tur\'an
function $\ex(n,\C F)$} is the maximum number of edges that an $\C F$-free
$k$-graph on $n$ vertices can have. This is one of the central questions of
extremal combinatorics that goes back to the fundamental
paper of Tur\'an~\cite{turan:41}. We refer the reader to the
surveys of the Tur\'an function by
F\"uredi~\cite{furedi:91}, Keevash~\cite{keevash:11},
and Sidorenko~\cite{sidorenko:95}.

As it was observed by Katona, Nemetz, and
Simonovits~\cite{katona+nemetz+simonovits:64}, 
the ratio 
$\ex(n,\C F)/{n\choose k}$ is non-increasing in $n$. In particular,
the limit
 $$
 \pi(\C F):=\lim_{n\to\infty} \frac{\ex(n,\C F)}{{n\choose k}}
 $$
 exists. It is called the \emph{Tur\'an density} of $\C F$. Let $\InfPi{k}$
consist of all possible Tur\'an densities of $k$-graph families and let
 $\FinPi{k}$ be the set of all possible Tur\'an densities when \textit{finitely many}
$k$-graphs are forbidden. Clearly, $\FinPi{k}\subseteq \InfPi{k}$.

For $k=2$, the celebrated Erd\H os-Stone-Simonovits 
Theorem~\cite{erdos+simonovits:66,erdos+stone:46} determines 
the Tur\'an density for every family $\C F$. In particular, we have
 \beq{PI2}
 \FinPi{2}=\InfPi{2}=\left\{\frac{m-1}m\mid m=1,2,3,\dots,\infty\right\}.
 \eeq
 (It is convenient to allow empty forbidden families, so $1\in \FinPi{k}$
for every $k$.)

Unfortunately, the Tur\'an function for
\emph{hypergraphs} (that is, $k$-graphs with 
$k\ge 3$) is much more difficult and many problems
(even rather basic ones) are wide open. 

Arguably, the case
when $|\C F|=1$ is the most interesting one. However, even
very simple forbidden hypergraphs turned out to be notoriously 
difficult. For example, the famous conjecture of Tur\'an
from 1941 that $\pi(\{K_4^3\})=5/9$ is still open, where $K_m^k$
denotes the complete $k$-graph on $m$ vertices. For no $3\le k<m$
is the value of $\pi(\{K_m^k\})$ known, despite the \$1000 prize of Erd\H os.
Razborov~\cite[Page~247]{razborov:11:psim} writes that \textit{``these questions
became
notoriously known ever since as some of the most difficult
open problems in discrete mathematics''}.

On the other hand, some Tur\'an-type results stop being true if only one
subgraph is to be forbidden. One such example is the Ruzsa-Szemer\'edi
theorem~\cite{ruzsa+szemeredi:78} that $\ex(n,\C F)=o(n^2)$, where
$\C F$ consists of all
$3$-graphs with $6$ vertices and at least $3$ edges (while $\ex(n,\{F\})=\Omega(n^2)$
for every $F\in\C F$). Some other problems 
(such as various intersection questions for uniform set systems, see e.g.~\cite{furedi:91}) can be restated in terms of the Tur\'an
function and require that more than one subgraph is forbidden. 
Also, new interesting
phenomena (such as, for example, non-principality, see~\cite{balogh:02,mubayi+pikhurko:08})
appear when one allows more than one forbidden $k$-graph. Last but not
least, by solving 
(perhaps more tractable) cases with $|\C F|>1$ we may get more insight
about the case $|\C F|=1$. In fact, some proofs that 
determine $\pi(\{F\})$ proceed by forbidding some extra hypergraphs
whose addition does not affect the Tur\'an density, see 
e.g.~\cite{baber+talbot:12,bollobas:74,frankl+furedi:89,mubayi:06,sidorenko:87}.

Little is known about $\FinPi{k}$ and $\InfPi{k}$ 
for $k\ge 3$. Brown and Simonovits~\cite[Theorem~1]{brown+simonovits:84}
noted that for every $\C F$ and $\e>0$ there is a 
finite $\C F'\subseteq \C F$ with $\pi(\C F')\le \pi(\C F)+\e$
(while, trivially, $\pi(\C F)\le \pi(\C F')$). It follows
that $\InfPi{k}$ lies in the closure of $\FinPi{k}$.
Here we show the following results about $\InfPi{k}$ with the first
one implying  that in fact $\InfPi{k}$ is the closure of $\FinPi{k}$.

\begin{proposition}\label{pr:closed} For every $k\ge 3$ the set 
$\InfPi{k}\subseteq [0,1]$ 
is closed. 
\end{proposition}

\begin{theorem}\label{th:uncountable}  For every $k\ge 3$ the set 
$\InfPi{k}$ has cardinality of the continuum.\end{theorem}

Since the number of finite families of $k$-graphs (up to isomorphism) 
is countable, Theorem~\ref{th:uncountable} implies that 
$\FinPi{k}\not=\InfPi{k}$ for $k\ge 3$, answering one part of a question
of Baber and Talbot~\cite[Question~31]{baber+talbot:12}.

Erd\H os~\cite{erdos:64} proved
that $\InfPi{k}\cap (0,k!/k^k)=\emptyset$,
that is, if the Tur\'an
density is positive, then it is at least $k!/k^k$. Let us call
a real $\alpha\in [0,1]$ a \emph{jump for $k$-graphs} 
if there is $\e>0$ such that $\InfPi{k}\cap (\alpha,\alpha+\e)=\emptyset$.
For example, every $\alpha\in [0,1]$ is a jump for graphs by~\req{PI2}
and every $\alpha\in [0,k!/k^k)$ is a jump for $k$-graphs by~\cite{erdos:64}.
The break-through paper of Frankl and R\"odl~\cite{frankl+rodl:84} 
showed that non-jumps exist for every $k\ge 3$, disproving the 
\$1000 conjecture of Erd\H os
 that $\InfPi{k}$ is
well-ordered with respect to the usual order on the reals. 
Further results on (non-)\,jumps were obtained
in~\cite{baber+talbot:11,frankl+peng+rodl+talbot:07,peng+zhao:08} and many other
papers.
 Our Theorem~\ref{th:uncountable} shows that $\InfPi{k}$ is
``very far'' from being well-ordered for $k\ge 3$. Since each jump is followed
by an interval disjoint from $\InfPi{k}$, at most
countably many elements of $\InfPi{k}$ can be jumps. Thus, by 
Theorem~\ref{th:uncountable},
the set of non-jumps has cardinality of the continuum.

Very few explicit numbers were proved to belong to $\FinPi{k}$. 
For example, before 2006 the only known members of $\FinPi{3}$ were $0$, 
$2/9$, $4/9$, $3/4$, and $1$ (see~\cite{bollobas:74,decaen+furedi:00,furedi+pikhurko+simonovits:03:ejc}). 
Then Mubayi~\cite{mubayi:06} showed that $(m-1)(m-2)/m^2\in \FinPi{3}$
for every $m\ge 4$. Very recently,
Baber and Talbot~\cite{baber+talbot:12} and Falgas-Ravry and 
Vaughan~\cite{falgas+vaughan:13}
determined
a few further elements of $\FinPi{3}$; their proofs are computer-generated,
being based on the flag algebra approach
of Razborov~\cite{razborov:10}. 
In all the cases when an explicit element of $\FinPi{k}$ is known, this 
limit density is achieved,
informally speaking,
by taking a finite pattern and blowing it up optimally.  
Here we generalise these results (as far as $\FinPi{k}$
is concerned) by showing that \emph{every}
finite pattern where, moreover, we are allowed to iterate the whole 
construction recursively inside a specified set of parts, produces 
an element of $\FinPi{k}$. 

Let us give some formal definitions. (We refer the reader to
Section~\ref{example} for an illustrative example.) A \emph{pattern}
is a triple $P=(m,E,R)$ where $m$ is a positive integer, $E$ is a
collection of $k$-multisets on $[m]:=\{1,\dots,m\}$, 
and $R$ is a subset of $[m]$. 
(By a \emph{$k$-multiset} we mean an unordered collection
of $k$ elements with repetitions allowed.) Let $V_1,\dots,V_m$ 
be disjoint sets and let $V=V_1\cup\dots\cup V_m$. The \emph{profile}
of a $k$-set $X\subseteq V$ (with respect to $V_1,\dots,V_m$) is
the $k$-multiset on $[m]$ that contains $i\in [m]$
with multiplicity $|X\cap V_i|$. For a $k$-multiset $Y\subseteq [m]$ let 
$\blow{Y}{V_1,\dots,V_m}$ consist of all $k$-subsets of $V$ whose profile is
$Y$. We call this $k$-graph the \emph{blow-up of $Y$} and
the $k$-graph
 $$
 \blow{E}{V_1,\dots,V_m}:=\bigcup_{Y\in E} \blow{Y}{V_1,\dots,V_m}
 $$ 
 is called the \emph{blow-up of $E$} (with respect to $V_1,\dots,V_m$).

A \emph{$P$-construction on a set $V$} is any $k$-graph $G$ that can
be recursively obtained as follows. Either let $G$ be the empty $k$-graph on $V$
(and
stop) or take an arbitrary 
partition $V=V_1\cup\dots\cup V_m$ where we require 
that if $i\in R$ then $V_i\not =V$.
 Add all edges of $\blow{E}{V_1,\dots,V_m}$ to $G$. Furthermore,
for every $i\in R$ take an arbitrary $P$-construction on $V_i$
and add all these edges to $G$. (If $R=\emptyset$, then there is
nothing to add and we have $G=\blow{E}{V_1,\dots,V_m}$.)
 Let $p_n$ be the maximum number of edges that can be obtained on $n$ vertices in this way:
 \beq{pn}
 p_n:=\max\big\{\,|G|\mid \mbox{$G$ is a $P$-construction on $[n]$}\,\big\}.
 \eeq

It is not hard to show (see Lemma~\ref{lm:lim})
that the ratio $p_n/{n\choose k}$ is non-increasing and therefore 
tends to a limit which
we denote by $\density{P}$ and call the \emph{Lagrangian}
of $P$:
 \beq{LagrangianP}
 \density{P}:=\lim_{n\to\infty} \frac{p_n}{{n\choose k}}.
 \eeq

 For $i\in [m]$ let $P-i$ be the pattern obtained
from $P$ by \emph{removing index $i$}, that is, we remove $i$ from $R$ and delete
all multisets containing $i$ from $E$ (and relabel the
remaining indices to form the set $[m-1]$). In other words, 
$(P-i)$-constructions are precisely those $P$-constructions where we
always let the $i$-th part be empty.  Let us call $P$ \emph{minimal}
if $\density{{P-i}}$ is strictly smaller than $\density{P}$ for every
$i\in [m]$. For example, the 2-graph pattern
$P:=(3,\{\,\{1,2\},\{1,3\}\,\},\emptyset)$ is not minimal as
$\Lambda_P=\Lambda_{P-3}=1/2$. 

\begin{theorem}\label{th:Fin} For every minimal pattern $P$ there is
a finite family $\C F$ of $k$-graphs such that for all $n\ge 1$ we
have $\ex(n,\C F)=p_n$ and, moreover, every maximum $\C F$-free $k$-graph on $[n]$
is a $P$-construction.\end{theorem}

\begin{corollary}\label{cr:Fin} For every pattern $P$ we have that 
$\density{P}\in\FinPi{k}$.\qed\end{corollary}

Corollary~\ref{cr:Fin} answers questions posed by Baber and 
Talbot~\cite[Question~29]{baber+talbot:12} and by 
Falgas-Ravry and Vaughan~\cite[Question~4.4]{falgas+vaughan:13}; we refer
the 
reader to Section~\ref{conclusion} for details.

Chung and Graham~\cite[Page~95]{chung+graham:eg} conjectured that $\FinPi{k}$
consists of rational numbers only. The following theorem 
disproves this conjecture for every $k\ge 3$. (Note that 
the conjecture is
true for $k=2$ by~\req{PI2}.) 
Independently, 
Chung and Graham's conjecture was disproved by Baber and 
Talbot~\cite{baber+talbot:12} who discovered a family 
of only three forbidden $3$-graphs whose Tur\'an density is irrational.
We should mention that Theorems~\ref{th:Fin} and~\ref{th:irrat} 
rely on the  Strong Removal Lemma of R\"odl and
Schacht~\cite{rodl+schacht:09} so they give 
families $\C F$ of
\textbf{huge} size.

\begin{theorem}\label{th:irrat} For every $k\ge 3$ the set
$\FinPi{k}$ contains an irrational number.\end{theorem}

This paper is organised as follows. Some further 
notation is given in Section~\ref{notation}.
The proof of Theorem~\ref{th:Fin} 
is presented in Section~\ref{Fin}; it is preceded by a number
of auxiliary results. Sections~\ref{irrat}, \ref{closed},
and~\ref{uncountable} contain the proofs of respectively Theorem~\ref{th:irrat},
Proposition~\ref{pr:closed}, and Theorem~\ref{th:uncountable}.
Finally, Section~\ref{conclusion} presents some concluding
remarks and open questions.

\section{Notation}\label{notation}

Let us introduce some further notation complementing and expanding 
that from the Introduction. Some other (infrequently used) definitions
are given shortly before they are needed for the first time in this
paper.

Recall that a $k$-\emph{multiset} $D$ is an unordered
collection of $k$ elements $x_1,\dots,x_k$ with repetitions
allowed. Let us denote this as $D=\multiset{x_1,\dots,x_k}$. The \emph{multiplicity $D(x)$} of $x$ in $D$ is the number of
times that $x$ appears. If the underlying set is understood to be 
$[m]$, then we can represent $D$ as the ordered
$m$-tuple $(D(1),\dots,D(m))$ of multiplicities. Thus, for example, the profile 
of $X\subseteq V_1\cup\dots\cup V_m$ is the multiset on $[m]$ whose
multiplicities are $(|X\cap V_1|,\dots,|X\cap V_m|)$. Also, let
$\rep{x}{r}$ denote the sequence consisting of $r$ copies of $x$; thus
the multiset consisting of $r$ copies of $x$ is denoted by $\multiset{\rep{x}{r}}$.
 If we need to
emphasise that a multiset is in fact a set (that is, no element
has multiplicity more than 1), we call it a \emph{simple set}.

For $D\subseteq [m]$ and sets $U_1,\dots,U_m$, denote $U_D:=\cup_{i\in D} U_i$.
Let ${X\choose m}:=\{Y\subseteq X\mid
|Y|=m\}$ consist of all $m$-subsets of a set $X$.
The \emph{standard $(m-1)$-dimensional simplex} is
 \beq{Sm}
 \IS_m:=\{\B x\in \I R^m\mid x_1+\dots+x_m=1,\ \forall\, i\in[m]\ x_i\ge 0\}.
 \eeq

\subsection{Hypergraphs}

We usually identify a $k$-graph $G$ with its edge set. For example,
$X\in G$ means that $X$ is an edge of $G$ and $|G|$ denotes
the number of edges. When we need to refer to the vertex set, we
write $V(G)$ and denote $v(G):=|V(G)|$. The \emph{(edge) density} of
$G$ is 
 $$
 \rho(G):=\frac{|G|}{{v(G)\choose k}}.
 $$ The \emph{complement of $G$} 
is $\OO G:=\{X\subseteq V(G)\mid |X|=k,\ X\not\in G\}$. For $x\in V(G)$ its
\emph{link} is the $(k-1)$-hypergraph
 $$
 G_x:=\{X\subseteq V(G)\mid x\not\in X,\ X\cup\{x\}\in G\}.
 $$
 For $U\subseteq V(G)$
its \emph{induced subgraph} is $G[U]:=\{X\in G\mid X\subseteq U\}$.
The vertex sets of $\OO G$, $G_x$, and $G[U]$ are by default $V(G)$,
$V(G)\setminus\{x\}$, and $U$ respectively. 
The \emph{degree} of $x\in V(G)$ is $d_G(x):=|G_x|$. Let $\Delta(G)$
and $\delta(G)$ denote respectively the maximum and minimum
degrees of the $k$-graph $G$.

An \emph{embedding} of a
$k$-graph $F$ into $G$ is an injection $f:V(F)\to V(G)$ such that
$f(X)\in G$ for every $X\in F$. An embedding is \emph{induced} if 
non-edges are mapped to non-edges.

\subsection{Pattern Specific Definitions}

Let $P=(m,E,R)$ be a pattern and $G$ be a $P$-construction on $[n]$. The
initial partition $V(G)=V_1\cup\dots\cup V_m$ is called the
\emph{level-$1$ partition} and $V_i$'s are called
\emph{level-$1$ parts}. For each $i\in R$
we denote the corresponding partition of $V_i$ as
$V_{i,1}\cup \dots\cup V_{i,m}$ and call these parts 
\emph{level-$2$ parts}. This notation generalises in the obvious 
way with $V_{i_1,\dots,i_s}$ for 
$(i_1,\dots,i_s)\in R^{s-1}\times [m]$ consisting of 
those vertices of $G$ that, for every $j=1,\dots,s$, 
belong to the $i_j$-th part on level~$j$. Also, we denote
$V_\emptyset:=V(G)$.

The \emph{length} of a sequence $\B i=(i_1,\dots,i_s)$ is $|\B i|=s$. The
sequence $\B i$ is \emph{legal} if $i_j\in R$
for all $j\in [s-1]$ and $i_s\in [m]$; this includes the empty 
sequence. 

We collect all parts that appear in the $P$-construction $G$ into 
a single vector
 $$
 \B V:=(V_\emptyset,V_1,\dots,V_m,\dots)
 $$ 
 and call $\B V$ the \emph{partition structure} of $G$; its index set 
is some subset of legal
sequences.

For convenience,
we view the  partition structure as vertical with
a level's index (called \emph{height}) 
increasing as we go up. In particular, 
the partition $V_1\cup\dots\cup V_m$ is called \emph{bottom}. 
By default, the profile of $X\subseteq V(G)$
is taken with respect to the bottom parts, that is, 
its multiplicities are $(|X\cap V_1|,\dots,|X\cap V_m|)$.
The \emph{branch $\branch{\B V}{x}$} of a vertex $x\in V(G)$ is the (unique) 
maximal sequence $\B i$ such that $x\in V_{\B i}$.

Given $P$, let $\C F_\infty$ consist of those $k$-graphs $F$ that do not embed into a
$P$-construction:
 \beq{CFInfty}
 \C F_\infty:=\{\mbox{$k$-graph $F$}\mid  
\mbox{every $P$-construction $G$ is $F$-free}\}.
 \eeq
 For an integer $n$, let $\C F_n$ consist of all members 
of $\C F_\infty$ with at most $n$ vertices:
 \beq{CFn}
 \C F_n:=\{F\in\C F_\infty\mid v(G)\le n\}.
 \eeq

Let the \emph{Lagrange polynomial} of $E$ be 
 \beq{LP}
 \lambda_E(x_1,\dots,x_m):=k!\,\sum_{D\in E}\; \prod_{i=1}^m\;
\frac{x_i^{D(i)}}{D(i)!}. 
 \eeq
 This definition is motivated by the fact that, 
for every partition $[n]=V_1\cup\dots\cup V_m$ we have that
 \beq{lambdaE}
 \rho(\blow{E}{V_1,\dots,V_m})
= \lambda_E\left(\frac{|V_1|}n,\dots,\frac{|V_m|}n\right)
+o(1),\qquad  \mbox{as $n\to\infty$};
 \eeq
 see also Lemma~\ref{lm:Omega1} that
relates $\lambda_E$ and $\density{P}$. 
The special case of~\req{LP} when $E$ is a $k$-graph (i.e.\ 
$E$ consists of simple sets)
has been successfully applied to Tur\'an-type problems, with 
the basic idea going back to
Motzkin and Straus~\cite{motzkin+straus:65}. Also, our 
definition of $\density{P}$ is a generalisation 
of the well-known \emph{hypergraph Lagrangian} 
$\Lambda_E:=\density{(m,E,\emptyset)}$,
see e.g.~\cite{baber+talbot:12}. 

For $i\in [m]$ let the \emph{link} $E_i$ consist of all $(k-1)$-multisets $A$
such that if we increase the multiplicity of $i$ in $A$ by one, then the
obtained
$k$-multiset belongs to $E$.
We call a pattern $P$ \emph{proper} if it is minimal and $0<\density{P}<1$.
Trivially,
every minimal pattern $P=(m,E,R)$ satisfies that 
 \beq{MinLinks}
 E_i\not=\emptyset,\qquad \mbox{for every $i\in [m]$.}
 \eeq

\subsection{An Example}\label{example}

To illustrate the above definitions, let us consider a specific simple example:
 \beq{example}
 P:=\big(\,2,\{\, \multiset{1,2,2}\, \}, \{1\}\,\big).
 \eeq
 Here a $P$-construction on $V$ is obtained by partitioning
$V=V_1\cup V_2$ with $V_1\not=V$ and adding all triples that have exactly two vertices
in $V_2$. Next, we apply recursion to $V_1$: namely, we partition $V_1=V_{1,1}\cup V_{1,2}$ 
with $V_{1,1}\not=V_1$ and add
all triples that intersect $V_{1,1}$ and $V_{1,2}$ in respectively one and two
vertices. Next, we repeat inside $V_{1,1}$, and so on. We can always stop; for example, we may choose to do this
after three iterations by letting $V_{1,1,1}$ span the empty 3-graph. 
In this case, the partition structure is 
 $$
 \B V=(V_\emptyset,V_1,V_2,V_{1,1},V_{1,2},V_{1,1,1},V_{1,1,2}),
 $$
 where $V_\emptyset:=V$.
 If we take a vertex $x$ in 
$V_2$, $V_{1,2}$, $V_{1,1,1}$, and $V_{1,1,2}$ then its branch is
respectively $(2)$, $(1,2)$, $(1,1,1)$, and $(1,1,2)$. This defines
the branch of every vertex as these four sets partition $V$; there is
no vertex whose branch is, for example, $(1,1)$. 

We have $\lambda_{E}(x_1,x_2)=6\cdot x_1\cdot (x_2^2/2)=3x_1x_2^2$. 
It is not hard to show (cf Lemma~\ref{lm:Omega1}) that $\Lambda_P=2\sqrt3-3$
and an example of a $P$-construction attaining this density is to use a ratio
close to $1:\sqrt3$ for each partition. Thus, this is an example of a 
3-graph pattern whose Lagrangian is irrational.

\section{Proof of Theorem~\ref{th:Fin}}\label{Fin}

The proof of Theorem~\ref{th:Fin} is rather long and relies on a number
of auxiliary results. 
Very briefly, it proceeds as follows. 
The starting point is the easy observation
(Lemma~\ref{lm:FInfty}) that by forbidding $\C F_\infty$ 
we restrict ourselves to $k$-graphs that embed into a $P$-construction;
thus $\ex(n,\C F_\infty)=p_n$. 
The deep and powerful Strong Removal Lemma of R\"odl and 
Schacht~\cite{rodl+schacht:09} (stated as Lemma~\ref{lm:RS} here) 
implies that for every $\e>0$
there is $M$ such that every $\C F_M$-free $k$-graph with
$n\ge M$ vertices can be made $\C F_\infty$-free by removing at most 
$\e {n\choose k}$ edges. It follows that
every maximum $\C F_M$-free $k$-graph $G$ on $[n]$ is $2\e {n\choose k}$-close 
in the edit distance to a $P$-construction,
see Lemma~\ref{lm:edit}. Although
the obtained $\e>0$ can be made arbitrarily small by choosing 
$M$ large, the author did not see any simple way of
ensuring that $\e\to 0$ for some fixed $M$ as $n\to\infty$.
Nonetheless our key Lemma~\ref{lm:Max} shows that some small but \emph{constant}
$\e>0$ suffices to ensure that there is a partition 
$V(G)=V_1\cup \dots\cup V_m$ such that 
$G\setminus (\cup_{i\in R} G[V_i])=\blow{E}{V_1,\dots,V_m}$, that is,
$G$ follows \emph{exactly} the bottom level of some $P$-construction 
(but nothing is stipulated about what happens inside the ``recursive''
parts $V_i$). 
The maximality of $G$ implies that each $G[V_i]$ with $i\in R$ is maximum
$\C F_M$-free (cf Lemma~\ref{lm:FmFree}), allowing us to
apply induction.

\subsection{Basic Properties of Patterns}\label{basic}

Here, let $P=(m,E,R)$ be an arbitrary pattern and let all
definitions of Sections~\ref{Intro} and~\ref{notation} apply.
In particular,  $p_n$, $\density{P}$, $\C F_\infty$
and $\C F_n$ are defined by respectively~\req{pn}, \req{LagrangianP},
\req{CFInfty}, and~\req{CFn}.

\begin{lemma}\label{lm:constr} Any induced subgraph (resp.\ any blow-up) of a
$P$-construction $H$
is (resp.\ embeds into) a $P$-construction.
\end{lemma}
 \bpf
 Let $\B V$ be the partition structure of $H$. If $H':=H[X]$ is an induced
subgraph, then we can initially
let $V'_{\B i}:=V_{\B i}\cap X$ for each index $\B i$. This need not
be a partition structure as we may have
$V_{i_1,\dots,i_s}'=V_{i_1,\dots,i_{s-1}}'$ for some $(i_1,\dots,i_s)\in R^s$, 
which is not allowed by the definition: namely, the partition
of $V_{i_1,\dots,i_{s-1}}'$ has the $i_s$-th part equal to the
whole set (and $i_s\in R$). We can fix one such occurrence by
removing $i_s$ from all indices that begin with
$(i_1,\dots,i_s)$. 
Formally, we remove all
parts $V_{i_i,\dots,i_{s-1},j_s,\dots}'$ with $j_s\not=i_s$ (note that they are all empty)
and relabel each part
$V_{i_1,\dots,i_s,j_{s+1},\dots,j_t}'$ into
$V_{i_1,\dots,i_{s-1},j_{s+1},\dots,j_t}'$.
We keep fixing all such occurences one by one. Since, 
for example, $\sum_{V_{\B i}'\in\B V'} |\B i|$ strictly
decreases each time, this procedure stops. The final
vector $\B V'$ shows that $H'$ is a $P$-construction.

If we insert a new vertex into a $P$-construction
by putting it into the same part as some existing vertex $x$, 
then we add all those edges (and possibly some further ones)
as when we just clone $x$. Thus every
blow-up of $H$,
which can be obtained by a sequence of
cloning steps and vertex removals, embeds into a $P$-construction.\epf

\begin{lemma}\label{lm:FInfty}  The following are equivalent for an arbitrary $k$-graph $G$ on $n$ vertices:
1) $G$ is $\C F_n$-free; 2) $G$ is $\C F_\infty$-free; 3) $G$ embeds into a $P$-construction; 4) $G$ embeds into a $P$-construction $H$ with $v(H)=n$.
\end{lemma}

\bpf The equivalence of 1), 2), and 3) follows from the definitions of $\C F_n$ and  $\C F_\infty$. Statements 3) and 4) are equivalent by Lemma~\ref{lm:constr}.\epf

It follows from Lemma~\ref{lm:FInfty} that $\ex(n,\C F_n)=\ex(n,\C F_\infty)=p_n$.

\begin{lemma}\label{lm:BlowInv} Let $s\in\I N\cup\{\infty\}$.
If $G$ is $\C F_s$-free, then any blow-up of $G$ is $\C F_s$-free.
\end{lemma}

\bpf Let $G'$ be obtained from $G$ by adding a clone $x'$ of some vertex
$x$ of $G$. Take any $U\subseteq V(G')$ with $|U|\le s$. If at least one
of $x$ and $x'$ is not in $U$, then $G'[U]$
is isomorphic to a subgraph of $G$ and cannot be in $\C F_s$; so suppose
otherwise. Since $G$ is $\C F_s$-free, there is an embedding $f$ of
$G[U\setminus\{x'\}]$ 
into some $P$-construction. By Lemma~\ref{lm:constr},
$G[U]$ is also embeddable.
It follows that $G'$ is
$\C F_s$-free.\epf

\begin{lemma}\label{lm:FmFree} Let $s\in\I N\cup\{\infty\}$. Let $G$ be a
$k$-graph on 
$V=V_1\cup\dots\cup V_m$ obtained by taking $\blow{E}{V_1,\dots,V_m}$
and putting arbitrary $\C F_s$-free $k$-graphs into parts $V_i$ with
$i\in R$. Then $G$ is $\C F_s$-free.\end{lemma}
 \bpf Take an arbitrary $U\subseteq V(G)$ with $|U|\le s$.
Let $U_i:=V_i\cap U$.
Note that $G[U_i]$ has no edges for $i\in [m]\setminus R$
and embeds into some $P$-construction
$H_i$ for $i\in R$ (because
$|U_i|\le s$ and $G[U_i]\subseteq G[V_i]$ is $\C F_s$-free). 
By combining
the partition structure of each $H_i$ together with
the level-1 decomposition $U=U_1\cup\dots\cup U_m$, we see
that $G[U]$ embeds into a $P$-construction, giving the required.\epf

\begin{lemma}\label{lm:lim} The ratio $p_n/{n\choose k}$
is non-increasing with $n$. In particular, the limit in \req{LagrangianP}
exists.
\end{lemma}

\bpf  Let $\ell< n$ and take a
maximum $P$-construction $G$ on $[n]$. Every $\ell$-subset
of $[n]$ spans at most $p_\ell$ edges by Lemma~\ref{lm:constr}. Averaging 
over all ${n\choose \ell}$ $\ell$-subsets gives that
$p_n\le p_\ell{n\choose \ell}/{n-k\choose \ell-k}=p_\ell{n\choose k}/{\ell\choose k}$,
as 
required.\epf

\begin{lemma}\label{lm:regular} For every $\e>0$ and $s\in\I N\cup\{\infty\}$
there is $n_0$
such that every maximum $\C F_s$-free 
$k$-graph $G$ with $n\ge n_0$ vertices has minimum degree at
least $(\density{P}-\e){n-1\choose k-1}$.
\end{lemma}
 \bpf Let $n$ be large and $G$ be as stated. Clearly, $|G|\ge p_n$.
The average degree of $G$ is $k|G|/n\ge kp_n/n\ge(\density{P}-\e/2){n-1\choose
k-1}$. If some $x$ has degree smaller than 
$(\density{P}-\e){n-1\choose k-1}$, then by deleting $x$ and adding a clone $y'$
of a vertex $y$ whose
degree is at least the average, we increase $|G|$ by at least
$|G_y|-|G_x|-{n-2\choose k-2}>0$. This preserves the $\C F_s$-freeness 
by Lemma~\ref{lm:BlowInv}, contradicting the
maximality of $G$.\epf

\begin{lemma}\label{lm:density=1}
 We have $\density{P}=1$ if and only if at least one of the following holds.
 \begin{enumerate}
 \item There is $i\in [m]$ such that $\multiset{\rep{i}{k}}\in E$;
 \item There are $i\in R$ and $j\in[m]\setminus\{i\}$ such that
$\multiset{\rep{i}{k-1},j}\in E$.
 \end{enumerate}
 \end{lemma}
 \bpf The converse implication is obvious: we can get the complete $k$-graph
on $[n]$ by taking $V_i=[n]$ in the first case and by taking $V_i=[n-1]$, $V_j=\{n\}$, and recursing inside $V_i$
in the second case.

Let us show the direct implication. Suppose that the above
multisets are not present in $E$. Let $n\to\infty$ and let $G$ be
a maximum $P$-construction on $[n]$ with the bottom partition $[n]=V_1\cup\dots\cup V_m$.

Suppose first that there is a part $V_i$ with $n-o(n)$ vertices 
for infinitely many $n$, say $i=1$. Assume that $1\in R$
for otherwise the complement $\OO G$ has at least 
${|V_1|\choose k}=\Omega(n^k)$ 
edges. Since $V_1$ is not allowed
to be the whole vertex set $[n]$, we can assume that e.g.\ $V_2\not=\emptyset$. 
Fix $x\in V_2$. The degree of $x$ in $G$ is at most 
$(n-|V_1|){n\choose k-2}=o(n^{k-1})$: since $\multiset{\rep{1}{k-1},2}\not 
\in E$, each edge of the link $(k-1)$-graph $G_x$ has to contain
at least one vertex outside of $V_1$. This contradicts
Lemma~\ref{lm:regular}.

Thus some two parts, say $V_1$ and $V_2$, have $\Omega(n)$
vertices each. Assume that $1\in R$ for otherwise at least
$\Omega(n^k)$ edges (those inside $V_1$) are missing from
$G$. Since $\multiset{\rep{1}{k-1},2}\not\in E$, all edges
that intersect $V_1$ in $k-1$ vertices and $V_2$ in one vertex
are not present. Again, at least $\Omega(n^k)$ edges are missing from $G$,
as required.\epf

The proof of Lemma~\ref{lm:density=1} shows that if $\density{P}=1$, then
the complete $k$-graph is a $P$-construction. This satisfies Theorem~\ref{th:Fin}
if we take $\C F=\emptyset$. Also, if $\density{P}=0$, then only empty 
$k$-graphs
are realisable as $P$-constructions and Theorem~\ref{th:Fin}
is also satisfied: let $\C F=\{K_k^k\}$ consist of a single edge. 
Thus it is enough to prove Theorem~\ref{th:Fin} for proper
patterns (that is, minimal patterns with Lagrangian strictly between 0 and 1).

\subsection{Properties of Proper Patterns}

In this section we let $P=(m,E,R)$ be an arbitrary pattern that is proper.
Here we establish some properties of $P$.

\begin{lemma}\label{lm:MaxVi} 
For every $P$-construction $G$ on $n$ vertices with
minimum degree $\delta(G)=\Omega(n^{k-1})$,
each bottom part $V_i$ has at most $(1-\Omega(1))n$ vertices
as $n\to\infty$.
\end{lemma}
\bpf For $i\in [m]\setminus R$ the claim follows from 
$\delta(G)n/k\le |G|\le {n\choose k}-{|V_i|\choose k}$.
Let $i\in R$. Since $V_i\not=V(G)$, pick
any vertex $x\in V_j$ with $j\not=i$. Since $\multiset{\rep{i}{k-1},j}
\not\in E$ by Lemma~\ref{lm:density=1}, 
every edge through $x$ contains at least one other
vertex outside of $V_i$. Thus $d_G(x)\le 
(n-|V_i|){n-2\choose k-2}$, implying the required.\epf

Let 
$$
 \IS_m^*:=\{\B x\in\I R^m\mid x_1+\dots+x_m=1,\ \forall\,i\in[m]\ 0\le x_i<1\}
 $$ 
 be obtained from $\IS_m$ by excluding the standard basis 
vectors, where $\IS_m$ is defined by~\req{Sm}.
Let us call a vector $\B x\in\I R^m$ \emph{optimal} if $\B x\in \IS_m^*$ and
 \beq{opt}
  \density{P}=\lambda_E(\B x) + \density{P} \sum_{i\in R} x_i^k.
 \eeq
 Let $\C X$ be the set of all optimal $\B x$. 
Note that when we define $\C X$ we restrict ourselves
to $\IS_m^*$ (i.e.\ we do not allow any standard basis vector to be included
into $\C X$). 

In a sense (with the formal statements appearing in
Lemmas~\ref{lm:Omega1} and~\ref{lm:OmegaN} below),
$\C X$ is precisely the set of optimal limiting ratios that
lead to asymptotically maximum $P$-constructions.
Let us illustrate this on the case when $P$ is
as in \req{example}. Suppose that we want to determine $\Lambda_P$.
Let $[n]=V_1\cup V_2$ be the bottom partition in a maximum $P$-construction $G$. Let 
$x_i:=|V_i|/n$ for $i=1,2$. By Lemma~\ref{lm:lim}, $\rho(G)$ and
$\rho(G[V_1])$ are close to $\Lambda_P$. (Note that
we cannot have $x_1=o(1)$ by Lemmas~\ref{lm:regular} and~\ref{lm:MaxVi}.) Thus, we conclude
that $\Lambda_P=3x_1x_2^2+\Lambda_P x_1^3+o(1)$, which is
exactly \req{opt} if we ignore the error term. Solving for $\Lambda_P$ and
excluding
$x_2$, we have to maximise $g(x):=3x(1-x)^2/(1-x^3)$ for $x\in (0,1)$.
In this particular case, the maximum is
$2\sqrt3-3$ and it is attained inside $(0,1)$ at the unique point 
$\alpha:=(\sqrt3-1)/2$.
It follows that $\Lambda_P=2\sqrt3-3$, \req{opt} has a unique solution in
$\IS_2^*$, and $\C X=\{(\alpha,1-\alpha)\}$. Note that although 
$(x_1,x_2)=(1,0)$ satisfies \req{opt}, we have that $\lim_{x\to 1}
g(x)=0<\Lambda_P$. This justifies why
we
exclude the standard basis vectors from $\C X$.

\begin{lemma}\label{lm:Omega1} Let $f(\B x):=\lambda_E(\B x) 
+ \density{P} \sum_{i\in R} x_i^k$ be
the right-hand side of~\req{opt}. Then the following claims hold.
 \begin{enumerate}
 \item\label{it:part1} $\C X\not=\emptyset$.
 \item\label{it:f} $f(\B x)\le \density{P}$ for
all $\B x\in \IS_m$. (Thus, by Part~\ref{it:part1}, $\C X$ is precisely the set of
elements in $\IS_m^*$ that maximise $f$.)
 \item\label{it:boundary} $\C X$ does not intersect the boundary 
of $\IS_m$.
 \item\label{it:f'} For every $\B x\in\C X$ and $j\in [m]$ we have
$\frac{\partial f}{\partial_j}(\B x)=k\density{P}$.
 \item\label{it:closed}  $\C X$ is a closed subset of $\IS_m$.
\item\label{it:almost} For every $\e>0$ there is $\alpha>0$ such that for every 
$\B y\in\IS_m$ with $\max(y_1,\dots,y_m)\le 1-\e$ and $f(\B y)\ge \density{P}-\alpha$ there
is $\B x\in\C X$ with $\|\B x-\B y\|_\infty\le \e$.
 \item\label{it:separated} There is $\beta>0$ such that for every 
$\B x\in \C X$ and every
$i\in[m]$ we have $x_i\ge \beta$.
 \end{enumerate}
 \end{lemma}
 \bpf Let $G$ be a maximum $P$-construction on $[n]$ with the bottom
partition $V_1\cup\dots\cup V_m$. By passing to
a subsequence of $n$, we can assume that, for every $i\in [m]$, the
ratio $|V_i|/n$ tends to some limit $x_i$. By 
Lemmas~\ref{lm:regular} and~\ref{lm:MaxVi}, 
$\B x=(x_1,\dots,x_m)$
belongs to $\IS_m^*$. By Lemma~\ref{lm:FmFree}, for each $i\in R$ the induced
subgraph $G[V_i]$ is a maximum $P$-construction. 
By Lemma~\ref{lm:lim}, we have that
$|G[V_i]|=(\density{P}x_i^k +o(1)){n\choose k}$. 
Now, \req{lambdaE} shows that $\B x$ satisfies~\req{opt}.
Thus $\B x\in \C X$, so this set 
is non-empty.

Let $\B x\in\IS_m$. If we use the
approximate ratios $x_1:\dots:x_m$ for the bottom partition 
$V_1\cup\dots\cup V_m$ and
put a maximum $P$-construction on each $V_i$ with $i\in R$, then
the obtained $P$-construction has edge density $f(\B x)+o(1)$. 
Thus $f(\B x)\le \density{P}$ for all $\B x\in \IS_m$, proving
Part~\ref{it:f}.

Suppose that $\C X$ intersects the boundary of $\IS_m$, that is,
$\C X$ contains some $\B x$
with zero entries. Without
loss of generality, assume that $x_1,\dots,x_{m'}$ are positive
while all other entries are 0.  Since $\B x\in\IS_m^*$, we have 
$m'\ge 2$.
Let $P'=(m',E',R')$ be obtained
from $P$ by removing the indices $m'+1,\dots,m$. 
Consider a $P'$-construction $H$ 
where the bottom partition $U_1\cup\dots\cup U_{m'}$ 
has approximate ratios $x_1:\dots:x_{m'}$ while each part
$U_i$ with $i\in R'$ spans a maximum $P'$-construction. By the
definition of $\density{P'}$, we have
 \begin{eqnarray*}
 \density{P'}&\ge& \rho(H)+o(1)\ =\ \lambda_{E'}(x_1,\dots,x_{m'})
+\density{P'}\sum_{i\in R'} x_i^k+o(1)\\
 &=& \lambda_{E}(\B x)
+\density{P'}\sum_{i\in R} x_i^k+o(1).
 \end{eqnarray*}
 Since $\B x\in\IS_m^*$, we have $\sum_{i\in R} x_i^k<1$. Thus
$\density{P'}+o(1)\ge \lambda_{E}(\B x)/(1-\sum_{i\in R} x_i^k)=\density{P}$,
 where we used identity~\req{opt} that $\B x\in\C X$
has to satisfy. This
contradicts the minimality of $P$ and proves Part~\ref{it:boundary}.

Let $\B x\in\C X$. By Part~\ref{it:boundary}, $\B x$ lies in 
the interior of $\IS_m$. Since $\B x$ maximises $f$ subject to
$x_1+\dots+x_m=1$, we conclude that all partial derivatives of 
$f$ coincide at $\B x$. Furthermore, this common value is $k\density{P}$,
which follows from the easy identity 
$\sum_{i=1}^m x_i\frac{\partial f}{\partial_i}(\B x)=kf(\B x)$,
establishing Part~\ref{it:f'}

From Part~\ref{it:f} we know that $\C X$  is precisely the
set of elements of $\IS_m^*$ that maximise $f(\B x)$.
Clearly, $f:\IS_m\to\I R$ is a continuous
function. Thus, in order
to prove Part~\ref{it:closed} it is enough to show that $\C X$
cannot accumulate to any element of the set 
$\IS_m\setminus\IS_m^*$ that consists of the standard basis vectors. 
The proof will essentially be a translation of the argument of Lemma~\ref{lm:MaxVi} into a more analytic language. Let $\B x\in\C X$.
Take any index $i\in [m]$. If $i\in [m]\setminus R$, then 
$\density{P}=f(\B x)\le 1-x_i^k+\Lambda_P(1-x_i)^k$, so $x_i$ cannot 
be arbitrarily close to $1$. Let
$i\in R$. 
Since $P$ is proper,  each 
monomial of $\lambda_E(\B x)$ contains at least two factors 
different from $x_i$ by Lemma~\ref{lm:density=1}.
Thus when we take the $j$-th derivative of $f$ for $j\in
[m]$, each monomial will have some factor
$x_s$ with $s\not=i$; of course, $x_s\le 1-x_i$. As the sum of
the coefficients of the degree-$k$ polynomial $f$ is, rather roughly, at most
$m^k$, we conclude
that $\frac{\partial f}{\partial_j}(\B x)\le km^k(1-x_i)$. By Part~\ref{it:f'} 
we conclude that
$1-x_i\ge \density{P}/m^k$, that is, $x_i$ is separated from $1$. This establishes
Part~\ref{it:closed}.

Suppose that Part~\ref{it:almost} is false. Then there is $\e>0$
such that for
every $i\in\I N$ there is $\B y_i\in\IS_m$ violating it with $\alpha=1/i$. By
the compactness of $\IS_m$ the sequence $(\B y_1,\B y_2,\dots)$
accumulates to some $\B y$. The vector $\B y$ belongs to $\IS_m^*$ 
by the assumption on each $\B y_i$. By the continuity of $f$ we have
$f(\B y)\ge \density{P}$, that is, $\B y\in\C X$, a contradiction
to $\B y$ being $\e$-far from $\C X$.

Part~\ref{it:separated} is proved in a similar way as Part~\ref{it:almost}.
(Alternatively, it directly follows from Parts~\ref{it:boundary}
and~\ref{it:closed}.)\epf

Informally speaking, the following lemma implies, among other things, 
that all part ratios 
of bounded height in a $P$-construction of large minimum degree
approximately follow some optimal vectors. For example, if $P$
is defined by~\req{example}, then Part~\ref{it:OmegaNx} of
Lemma~\ref{lm:OmegaN} gives that,
for any fixed $\ell$, every $P$-construction $G$ on $[n]$
with $\delta(G)\ge (\Lambda_P-o(1)){n-1\choose k-1}$ satisfies
$|V_{\rep{1}{s},2}|:|V_{\rep{1}{s},1}|=\sqrt3+o(1)$ for each
$0\le s\le\ell$.

\begin{lemma}\label{lm:OmegaN} 
For every $\e>0$ and $\ell\in\I N$ there are constants
$\alpha_0,\e_0,\dots,\alpha_\ell,\e_\ell,\alpha_{\ell+1}\in (0,\e)$ 
and $n_0\in\I N$ such that the following holds. Let $G$ be 
an arbitrary $P$-construction $G$ on $n\ge n_0$ vertices with the partition
structure $\B
V$
such that the minimal degree $\delta(G)\ge (\density{P}-\alpha_0){n-1\choose
k-1}$. Take arbitrary $\B i\in R^{s}$ with $0\le s\le \ell$ and denote $\B v_{\B
i}:=(|V_{\B i,1}|/|V_{\B i}|,\dots, |V_{\B
i,m}|/|V_{\B i}|)$. Then:
 \begin{enumerate}
 \item\label{it:OmegaNx}  $\|\B v_{\B i}-\B x\|_\infty\le \e_s$ for
some $\B x\in\C X$;
 \item \label{it:OmegaNVi} $|V_{\B i,j}|\ge (\beta/2)^{s+1}n$ for all $j\in[m]$,
where $\beta$ is returned by Part~\ref{it:separated}
of Lemma~\ref{lm:Omega1};
 \item\label{it:OmegaNDelta} $\delta(G[V_{\B i,j}])\ge 
(\density{P}-\alpha_{s+1}){|V_{\B i,j}|-1\choose k-1}$ for all $j\in R$.
 \end{enumerate}
\end{lemma}

\bpf We choose positive constants in this order 
 $$
 \alpha_{\ell+1}\gg \e_\ell\gg \alpha_\ell\gg 
\dots\gg \e_0\gg
\alpha_0\gg 1/n_0,
 $$
 each being sufficiently small depending on $P$, $\e$, $\beta$, and the previous constants. Let $G$ and $\B i$ be as in the lemma.
We use induction on $s=0,1,\dots,\ell$. Let $U_j:=V_{\B i,j}$ for $j\in
[m]$, $U:=V_{\B i}=U_1\cup\dots\cup U_m$, and $\B u:=\B v_{\B
i}=(|U_1|/|U|,\dots,|U_m|/|U|)$. By
the inductive assumption on $\delta(G[U])$ (or by the assumption
of the lemma if $s=0$ when $U=V_\emptyset=V(G)$), we have that
 $$
 |G[U]|\ge \delta(G[U])|U|/k\ge 
(\density{P}-\alpha_s){|U|\choose k}.
 $$ 
 Since $\alpha_s\ll \e_s$, Lemma~\ref{lm:MaxVi} and
Part~\ref{it:almost} of Lemma~\ref{lm:Omega1} (when applied to $\B y=\B u$) 
give the desired $\B x$, proving Part~\ref{it:OmegaNx}.

For all $j\in[m]$, we have $|U_{j}|\ge
(x_j-\e_s)|U|
\ge (\beta/2)|U|$, which is at least 
$(\beta/2)^{s+1}n$ by the inductive assumption, proving
Part~\ref{it:OmegaNVi}.

Finally, take arbitrary $j\in R$ and 
$y\in U_j$. The degree of $y$ in $\blow{E}{U_1,\dots,U_m}$ is
 $$
 d_{\blow{E}{U_1,\dots,U_m}}(y)=\left(\frac1k\times\frac{\partial
\lambda_E}{\partial_j}(\B u)+o(1)\right) {|U|-1\choose k-1}.
 $$ 
Since $\|\B u-\B x\|_\infty\le\e_s\ll \alpha_{s+1}$, we have by Part~\ref{it:f'}
of Lemma~\ref{lm:Omega1} that, for example,
 $$
 \frac{\partial}{\partial_j}\lambda_E(\B u)-\alpha_{s+1}^2\le 
\frac{\partial}{\partial_j}\lambda_E(\B x)=\frac{\partial}{\partial_j} f(\B
x)-k\density{P} x_j^{k-1}= k\density{P}-k\density{P} x_j^{k-1}.
 $$ 
 Thus, by the (inductive) assumption on the minimal degree of $G[U]$, we have
 \begin{eqnarray*}
 d_{G[U_j]}(y) &=&  d_{G[U]}(y)-d_{\blow{E}{U_1,\dots,U_m}}(y)\\
 &\ge& \left((\density{P}-\alpha_s)- (\density{P}-\density{P} x_j^{k-1}+
2\alpha_{s+1}^2/k)\right)
{|U|-1\choose k-1}\\
 &\ge & \left(\Lambda_P(u_j-\e_s)^{k-1} - 3\alpha_{s+1}^2/k\right){|U|-1\choose
k-1}\ \ge\ (\Lambda_P-\alpha_{s+1}){|U_j|-1\choose k-1},
 \end{eqnarray*}
 where we used $|u_j-x_j|\le\e_s$ and
$x_j\ge\beta\gg\alpha_{s+1}\gg \e_s\gg \alpha_s$. This finishes the proof of Lemma~\ref{lm:OmegaN}.\epf

Recall that the link $E_i$ of $i\in [m]$ consists of all $(k-1)$-multisets on $[m]$ such that if we increase the multplicity of $i$ by one, then the obtained $k$-multiset belongs to $E$.

\begin{lemma}\label{lm:DNs} If distinct $i,j\in[m]$ satisfy $E_i\subseteq E_j$, then
$i\in R$, $j\not\in R$, and $E_i\not=E_j$. In particular, no two vertices of the
pattern $P=(m,E,R)$ have the same links in $E$.
\end{lemma}
 \bpf Take some optimal $\B x\in \C X$. By
Part~\ref{it:boundary} of Lemma~\ref{lm:Omega1}, all coordinates
of $\B x$ are non-zero. Define $\B x'\in\IS_m$ by $x_i'=0$,
$x_j'=x_i+x_j$, and $x_h'=x_h$ for all other indices $h$. 
We claim that
 \beq{lambdaxx'}
 \lambda_E(\B x')\ge\lambda_E(\B x).
 \eeq One way to
show \req{lambdaxx'} is to use \req{lambdaE}. Consider
some $F:=\blow{E}{V_1,\dots,V_m}$.
The assumption $E_i\subseteq E_j$ implies
that if decrease the multiplicity of $i$ in some $A\in E$ but
increase the multiplicity of $j$ by the same amount, then the new multiset
necessarily belongs to $E$. Thus if we remove a vertex $y$ from $V_i$ and
add a vertex $y'$ to $V_j$, then the obtained $k$-graph $F'$ has at least
as many edges as $F$. (In fact, we have that $F_y\subseteq F_{y'}'$.)
 Since $\B x'$ is obtained from $\B x$ by
shifting weight from $x_i$ to $x_j$, \req{lambdaxx'} follows. 

Also, $m\ge 3$ for otherwise $\multiset{\rep{j}{k}}\in E$, contradicting Lemma~\ref{lm:density=1}. Thus $\B x'\in\IS_m^*$.

We conclude that $j\not\in R$ for otherwise we get a contradiction to
Part~\ref{it:f} of Lemma~\ref{lm:Omega1} by using \req{lambdaxx'}
and the trivial inequality $x_i^k+x_j^k<(x_i+x_j)^k$.
Likewise, $i\in R$ for otherwise the vector $\B x'$, that has a zero 
coordinate, would belong to $\C X$.
Finally, we see  that $E_i\not=E_j$ by swapping
the roles of $i$ and $j$ in the above argument.\epf

A map $h:[m]\to[m]$ is
an \emph{automorphism} of the pattern $P$ if $h$ is bijective, $h(R)=R$, and $h$ is
an automorphism of $E$ (that is, $h(E)=E$). Let us call a $P$-construction $G$
with the bottom partition $V_1\cup\dots\cup V_m$ \emph{rigid} if 
for every embedding $f$ of $G$ into a $P$-construction 
$H$ with the bottom
partition $U_1\cup\dots\cup U_m$ such that $f(V(G))$ intersects at least
two different parts $U_i$, there is an automorphism $h$ of $P$ 
such that $f(V_i)\subseteq U_{h(i)}$
for every $i\in[m]$. 

For example, the pattern $P$ in~\req{example}
has no non-trivial automorphism and a rigid $P$-construction can 
be obtained by taking any $\blow{E}{V_1,V_2}$ with $|V_1|\ge 1$ and $|V_2|\ge
3$. Thus the following Lemma~\ref{lm:MaxRigid} is trivially true for
this
particular $P$.

\begin{lemma}\label{lm:MaxRigid} For all large $n$, every
maximum  $P$-construction
$G$ on $[n]$ is rigid.
 \end{lemma}

Since the proof of Lemma~\ref{lm:MaxRigid} in general is long and complicated, 
some informal discussion may
be helpful here. 
It is not surprising
that the proof is far simpler if $R=\emptyset$. In fact, an 
example of a rigid $P$-construction in this case can be obtained by letting each $V_i$
have more than $(k-2)m$ vertices. Indeed, take any embedding $f$
of $G=\blow{E}{V_1,\dots,V_m}$ 
into $\blow{E}{U_1,\dots,U_m}$. For every $i\in [m]$ at least $k-1$
vertices of $V_i$ go into $U_{h(i)}$ for some $h(i)$. 
It is not hard to see that 
if we map each part $V_i$ entirely into $U_{h(i)}$,
then the new map is also an embedding. 
Since $P$ is minimal, $h$ has to be surjective and some
extra work shows that necessarily $f(V_i)\subseteq U_{h(i)}$
for all $i\in[m]$. 
(In fact, if furthermore $E$ consists of
simple $k$-sets only, then $|V_i|\ge 1$ is enough for rigidity.) 

The case $R\not=\emptyset$ is more complicated, although the main ideas
(such as using the function $h$ that specifies where a large part of 
$V_i$ is mapped to) are roughly the same. One complication is that
for a non-minimal pattern $P$ 
there can be embeddings that map the bottom edges into
different levels. For example, let 
 $$
 P=\big(5,\left\{\,
 \multiset{1,2,3},
 \multiset{1,2,4},
 \multiset{1,2,5},
 \multiset{3,4,5}\,\right\},\{5\}\big)
 $$
 and let $f$ map the bottom parts $V_1,\dots,V_5$ into respectively 
$U_1,U_2,U_{3,1},U_{3,2},U_{3,3}$. Here,  $P$
is obtained from the pattern $(3,\{\, \multiset{1,2,3}\, \}, \{3\})$
by ``expanding'' the third part up one level. Thus our proof of 
Lemma~\ref{lm:MaxRigid}	 should in particular catch all such redundancies.

\bpf[Proof of Lemma~\ref{lm:MaxRigid}.] Let $n\to\infty$
and $G$ be a maximum $P$-construction on $[n]$ with
the partition structure $\B V$. Take any 
embedding $f$ of $G$ into some $P$-construction 
$H$ with the bottom partition
$V(H)=U_1\cup\dots\cup U_m$ such that $f(V(G))$ intersects at least two
different parts $U_i$.

\bcl{induced}{The map $f$ is an induced embedding (that is, $f(X)$ is an edge if and
only if $X$ is).}

\bcpf If some non-edge $D\in \OO G$ is mapped by $f$ into an edge of $H$, then
the $k$-graph $G\cup\{D\}$ embeds into a $P$-construction (the
very same map $f$ embeds it into $H$). However, this   contradicts the maximality
of $G$.\ecpf

By Lemma~\ref{lm:regular} and Part~\ref{it:OmegaNVi} of Lemma~\ref{lm:OmegaN},
the size of 
each $V_i$ tends to infinity. By the pigeonhole principle,
there is a function $h:[m]\to [m]$ such that 
 \beq{h(i)}
 |f(V_i)\cap U_{h(i)}|\ge k,\quad \mbox{for all $i\in [m]$}.
 \eeq

\bcl{h}{We can choose $h$ in~\req{h(i)} so that,
additionally, $h(R)\subseteq R$ and $h$
assumes at least two different values.}

\bcpf Suppose that $R\not=\emptyset$ and 
we cannot satisfy the first part of the 
claim for some $i\in R$, that
is, for each $s\in R$ we have $|f(V_i)\cap
U_s|<k$. Thus $G[V_i]$ with exception
of at most $(k-1)|R|$ vertices is embeddable into $H[U_{[m]\setminus R}]$. 
By the maximality of $G$, Lemmas~\ref{lm:regular} and~\ref{lm:OmegaN}
give that $|V_i|\to\infty$ and $\rho(G[V_i])= \density{P}+o(1)$. 
This means that $(P-R)$-constructions can
contain arbitrarily large subgraphs of edge density 
$\density{P}+o(1)$, that is, $\density{P-R}\ge
\density{P}$. However, this contradicts the minimality of $P$.

Let us restrict ourselves to those $h$ with $h(R)\subseteq R$. 
Suppose that we cannot fulfil the second part of the claim. 
Then there is $j\in [m]$ such that
$|f(V_i)\cap U_j|\ge k$ for every $i\in[m]$. 
Since $E\not=\emptyset$, the induced subgraph
$G[f^{-1}(U_j)]$ is non-empty (it has at least $k$ vertices
from each $V_i$) and is mapped entirely into $U_j$. Thus
$j\in R$. Since $f(V(G))$
intersects at least two different parts $U_i$, we can pick some $x\in
V_i$ with $f(x)\in U_s$ and $s\not=j$.  Fix some $(k-1)$-multiset
$D\in E_i$. (Note that $E_i\not=\emptyset$ by \req{MinLinks}.)
Take an edge $D'\ni x$ of $G$ so that $D'\setminus\{x\}$
is a subset of $f^{-1}(U_j)$ and has profile $D$; it exists because
each part $V_g$ contains at least $k$ vertices of $f^{-1}(U_j)$.
The $k$-set $f(D')$ is an edge of $H$ as $f$ is an embedding. However, it has
$k-1$ vertices in $U_j$ and one vertex in $U_s$. Thus the
$k$-multiset $\multiset{\rep{j}{k-1},s}$ 
belongs to $E$. Since $j\in R$, this contradicts 
Lemma~\ref{lm:density=1}. The claim is proved.\ecpf

\bcl{bijection}{Each $h$ satisfying Claim~\ref{cl:h} is a bijection.}

\bcpf For $j\in [m]$ let $U_{j}':=\cup_{i\in h^{-1}(j)} f(V_i)\subseteq V(H)$. 
(Thus $U_j'=\emptyset$
for $j$ not in the image of $h$.) Let $H'$ be the $P$-construction
on $f(V(G))$ such that $U_1'\cup\dots\cup U_m'$ is the bottom partition of $H'$ 
and,  for $i\in R$, $H'[U_i']$ is the image of the $P$-construction 
$G[f^{-1}(U_i')]$ under the bijection $f$. 

We have just defined a new $P$-construction $H'$
so that each part $V_i$ of $G$ is entirely mapped by $f$ into the
$h(i)$-th part of $H'$, that is, \emph{all} vertices of $G$ follow $h$ now.
This $H'$ will be used only for proving  that $h$ is a bijection. 
The reader should be able to derive from the proof of Claim~\ref{cl:bijection} that in fact
$U_j'\subseteq U_j$ and $H'=H[f(V(G))]$ (but we will not use
these properties). 

Let us show first that the same map $f$
is an embedding of $G$ into $H'$.  First, take any bottom edge $D\in G$ such that $f(D)$ intersects two different
parts $U_i'$. Let $D'\in G$ have the same profile as $D$ and satisfy
 \beq{D'}
 D'\subseteq \cup_{i\in[m]} \left( V_i\cap f^{-1}(U_{h(i)}) \right),
 \eeq
 which is possible because there are at least $k$ vertices available in
each part $V_i$. Since $f(D')\cap U_i= f(D')\cap U_i'$ for all $i\in[m]$, 
the $f$-image of  $D'$ has the same profile $X$ with
respect to the partitions $U_1\cup\dots\cup U_m$ and
$U_1'\cup\dots\cup U_m'$. Thus $X\in E$.
Next, as each $f(V_i)$ lies entirely inside $U_{h(i)}'$, the
sets $f(D)$ and $f(D')$ have the same profiles with
respect to parts $U_i'$.
Thus $f(D)$ is an edge of $\blow{E}{U_1',\dots,U_n'}$, 
as required. Next, take any $i\in[m]$ and let $G':=G[f^{-1}(U_i')]$. Assume that
$i\in [m]\setminus R$ for otherwise $f(G')\subseteq H'$ by the
definition of $H'$. We claim that $G'$ has no edges in this case. Since
$h(R)\subseteq R$, we have
$h^{-1}(i)\cap R=\emptyset$. Thus it remains to derive a contradiction
by assuming that a bottom edge $D$ of $G$ belongs to $G'$. As before,
we can find an edge $D'\in G$ that satisfies~\req{D'} and has the same 
profile as $D'$ with respect to $V_1,\dots,V_m$.
However, $f$ maps this $D'$ inside a non-recursive part $U_i$ of $H$,
a contradiction. Thus $f$ is an embedding of $G$ into $H'$.

Thus, by considering $H'$ instead of $H$ (and without changing
$h$) we have that $f(V_i)\subseteq U_{h(i)}'$ for all $i\in[m]$. 

Suppose on the contrary to the claim that $|h^{-1}(s)|\ge 2$ for some $s\in [m]$.  Let $A:=h^{-1}(s)$
and $B:=[m]\setminus A$. Since $h$ assumes at least two different values,
the set $B$ is non-empty.

Note that
$U_s'$ is \emph{externally $H'$-homogeneous}, meaning that any
permutation $\sigma$ of $V(H')$ that fixes every vertex outside of $U_s'$ is a
symmetry of the set of $H'$-edges that intersect the complement of $U_s'$, that
is, $\sigma\big(H'\setminus {U_s'\choose k}\big)=H'\setminus {U_s'\choose k}$.
It follows from Claim~\ref{cl:induced} that $f^{-1}(U_s')=V_A$
is
externally $G$-homogeneous. 
(Recall that we denote $V_A:=\cup_{i\in A} V_i$.) 
Since each $V_i$ has at least $k$ elements, we conclude that $A$ is externally $E$-homogeneous.

Suppose first that $A\cap R\not=\emptyset$. By the above homogeneity, if we
replace $G[V_A]$ by any $P$-construction, then the new $k$-graph on $V$ is
still a $P$-construction. Also, recall that each $V_i$ has
size $\Omega(n)$ by Lemmas~\ref{lm:regular} and~\ref{lm:OmegaN}. 
Thus, by the maximality of
$G$, the edge density of $G[V_A]$ is $\density{P}+o(1)$. Also, 
$\rho(G[V_i])=\density{P}+o(1)$ for $i\in A\cap R$. Consider
the pattern $Q:=P-B$ obtained by removing $B$ from $P$. Without loss
of generality assume that
$A=[a]$. For $i\in A$
let $x_i:=|V_i|/|V_A|$. The obtained vector $\B x\in \IS_{a}$
satisfies
$\density{P}= \lambda_Q(\B x)+ \sum_{i\in A\cap R} \density{P}
x_i^k+o(1)$.  On the other hand, if we use the same vector $\B x$
for the bottom ratios and put a
maximum $Q$-construction on each recursive part, then we
get
overall density at most $\density{Q}+o(1)$. Thus $\density{Q}\ge \lambda_Q(\B x)+ 
\sum_{i\in A\cap R} \density{Q} x_i^k+o(1)$. Since $G$ is a maximum
$P$-construction, we have that each part $V_i$ has $\Omega(n)$ vertices; 
thus no $x_i$ can be equal $1-o(1)$ and we have that $1-\sum_{i\in A\cap R}
x_i^k=\Omega(1)$. These inequalities imply
that $\density{Q}\ge \lambda_Q(\B x)/(1-\sum_{i\in A\cap R}
x_i^k)+o(1)=\density{P}+o(1)$, contradicting the minimality of $P$.

Finally, suppose that $A\cap R=\emptyset$. Since $V_A$ is externally
$G$-homogeneous
and $A$ consists of at least two indices $i\not=j$, we have that $E$ contains at least one multiset entirely inside $A$ (for 
otherwise $E_i=E_j$, contradicting Lemma~\ref{lm:DNs}). Since
$f(V_A)= U_s'$, we have that $s\in R$. By the maximality of $G$
and Claim~\ref{cl:induced} it follows that
the edge density of $H'[U_s']$ (and thus of $G[V_A]$) is $\density{P}+o(1)$. Thus
$\density{P-B}\ge \density{P}$, a contradiction proving the claim.\ecpf

It follows from Claim~\ref{cl:bijection} that each $h$ satisfying Claim~\ref{cl:h} is an automorphism of $P$. By 
relabelling the parts of $H$, we can assume 
for notational convenience that $h$ is the identity mapping. 
Now we are ready to prove the lemma,
namely that $f(V_i)\subseteq U_i$ for every $i\in [m]$.

Suppose on the contrary that $f(x)\in U_j$ for some $x\in V_i$
and $j\in [m]\setminus\{i\}$. It follows 
that $E_i\subseteq E_j$. By Lemma~\ref{lm:DNs} this inclusion is
strict and $i\in R$. Pick some $X$ from $E_j\setminus E_i\not=\emptyset$. We can find
$D\in H$ such that $D\setminus\{f(x)\}$ has the profile $X$ with respect to
both $U_1\cup \dots\cup U_m$
and $f(V_1)\cup\dots\cup f(V_m)$.
But then $f^{-1}(D)$ is not an edge of $G$ because its profile $X\cup\{i\}$
is not in $E$. (Note that it is impossible that $X=\multiset{\rep{i}{k-1}}$ as
this
would give that $E$ contains $\multiset{\rep{i}{k-1},j}$, the profile of $D\in H$, 
contradicting Lemma~\ref{lm:density=1}.) 
Thus $f$ is not induced, contradicting Claim~\ref{cl:induced}. 
This shows that $G$ is rigid.\epf

\begin{lemma}\label{lm:AddRigid} Every rigid 
$P$-construction $G$ with the partition structure
$\B V$ such that $|V_{\B i}|\ge k$ for
every legal $\B i$ with $|\B i|\le 2$ remains rigid after
the addition of any new vertex.
\end{lemma}

\bpf It is enough to show that if we take any embedding
$f$ of $G$ into a $P$-construction with the partition structure $\B U$
such that $f(V_i)\subseteq U_i$ for $i\in[m]$ and 
add one vertex $x$ to some part $V_i$, then
any extension of $f$ to $x$ maps it into $U_i$.
If $i\in [m]\setminus R$, then $x$ and some $y\in V_i$ 
have the same links in $V':=V(G)\setminus \{x,y\}$; since the part containing 
$f(y)$ is determined by the values of $f$ on $V'$, the same applies
to $f(x)$, as required. So let $i\in R$. Since $|V_{i,j}|\ge k$ for each $j\in [m]$, 
the link of $x$ in $G$ necessarily contains at least one $(k-1)$-set entirely
inside $V_i$ by~\req{MinLinks}. This forces by Lemma~\ref{lm:density=1} that
$f(x)\in U_i$,
finishing the proof.\epf

Later (in the proof of Lemma~\ref{lm:Max}) we will need
the existence of a rigid $P$-construction
such that the recursion goes for exactly $\ell$ levels for some $\ell$ 
and every part at height at most $\ell$
has many vertices. This can be achieved as follows.
Take large $n$ and let $G$ be a maximum $P$-construction
on $[n]$. It is rigid by Lemma~\ref{lm:MaxRigid}. 
Also,
by Lemma~\ref{lm:regular} and Part~\ref{it:OmegaNVi} of Lemma~\ref{lm:OmegaN},
$G$ satisfies the assumptions of Lemma~\ref{lm:AddRigid}. 
Thus we can add some extra vertices to the $P$-construction $G$, without
increasing its maximum height $\ell$ while achieving the following property:

\begin{lemma}\label{lm:tight} There are $\ell\in\I N$ and a rigid
$P$-construction
with the partition structure $\B V=(V_{\B i}\mid |\B i|\le \ell)$ 
such that for every legal sequence $\B i$ of length at most $\ell$
we have $|V_{\B i}|\ge (k-1)\max(m,k)$.\qed\end{lemma}

\subsection{Key Lemmas}

In this section, $P=(m,E,R)$ is still a proper pattern.
Let us call two $k$-graphs with the same number of vertices \emph{$s$-close}
if one can be made isomorphic to the other by changing at most $s$ edges.

\begin{lemma}\label{lm:Max} There are $c_0>0$ and $M_0\in\I N$ such that
the following holds.
Let $G$ be a maximum $\C F_{M_0}$-free $k$-graph on $n\ge 2$ vertices that is
$c_0{n\choose k}$-close to some $P$-construction. Then there is a partition
$V(G)=V_1\cup\dots\cup V_m$ such that no $V_i$ is equal to $V$ and
 \beq{BottomLevel}
 G\setminus (\cup_{i\in R} G[V_i])=\blow{E}{V_1,\dots,V_m}.
 \eeq
\end{lemma}

\bpf Clearly, it is enough to establish the existence of $M_0$ 
such
that the conclusion of the lemma holds for every sufficiently large $n$. (Indeed,
it clearly holds for $n\le M_0$ by Lemma~\ref{lm:FInfty}, so we can simply increase $M_0$
at the end to take care of finitely many exceptions; alternatively, one 
can decrease $c_0$.)

Let $F$ be the rigid construction returned by Lemma~\ref{lm:tight}.
Let $\ell$ be the maximum height of $F$ and let $\B W$ be its partition structure.
(Our proof also works if $R=\emptyset$, when $\ell=1$; in fact,
some parts can be simplified in this case.)
Let $M_0:=v(F)+k$.

We choose some constants $c_i$ in this order $c_4\gg c_3\gg c_2\gg
c_1\gg c_0>0$, each being sufficiently small depending on the previous ones. Let
$n$ tend to infinity.

Let $G$ be a maximum $\C F_{M_0}$-free $k$-graph on $[n]$
that is $c_0{n\choose k}$-close to some $P$-construction $H$. We
can assume that $V(H)=[n]$ and the vertices of $H$ are already 
re-labelled so that $|G\bigtriangleup H|\le c_0{n\choose k}$. 
Let $\B V$ be the partition structure of $H$. 
In particular, the bottom partition of $H$ is $V_1\cup\dots\cup V_m$.

One of the technical difficulties that we are going to face is that
some part $V_i$ with $i\in R$ may in principle contain almost every vertex of $V(G)$ (so every other part $V_j$ has $o(n)$ vertices).
This means that the ``real'' approximation to $G$ starts only at some
higher level inside $V_i$. On the other hand, Lemma~\ref{lm:OmegaN} 
gives us a way to rule out such cases: we have to ensure that
the minimal degree of $H$ is close to $\density{P}{n-1\choose
k-1}$. So, as our first step, we are going to modify the 
$P$-construction $H$ (perhaps at the expense of increasing
$|G\bigtriangleup H|$ slightly) so that its minimal degree is large.

Namely, let $Z:=\{x\in [n]\mid d_H(x)< (\density{P}-2c_1) {n-1\choose k-1}\}$.
By Lemma~\ref{lm:regular} we can assume that 
 \beq{DeltaG}
 \delta(G)\ge (\density{P}-c_1) {n-1\choose k-1}.
 \eeq 
 Thus every vertex of $Z$ contributes at least 
$c_1{n-1\choose k-1}/k$ to $|G\bigtriangleup H|$. We conclude that
$|Z|\le c_0n/c_1$. Fix an arbitrary $y\in [n]\setminus Z$. 
Let us change $H$ by making all vertices in
$Z$ into clones of $y$ (and
updating $\B V$
accordingly). Clearly, we have now
 \beq{MinDegStage1}
 \delta(H)\ge (\density{P}-2c_1){n-1\choose k-1}-|Z|{n-2\choose k-2}\ge
 (\density{P}-3c_1){n-1\choose k-1}
 \eeq
 while $|G\bigtriangleup H|\le c_0{n\choose k}+ |Z|{n-1\choose k-1}\le
c_1{n\choose k}$.

If we end up with an improper partition structure (e.g.\ $V_i=V(H)$ for some $i\in R$),
then we correct this as in the proof of Lemma~\ref{lm:constr} without
changing the $k$-graph $H$. 

By Lemma~\ref{lm:OmegaN} we can conclude that (in the new $k$-graph $H$) all
part ratios up to
height $\ell$ are close to optimal ones and $|V_{\B i}|\ge 2c_4n$ for each legal sequence $\B i$ of length at most~$\ell$.

Let $A:=\blow{E}{V_1,\dots,V_m}\setminus G$
consist of what we shall call \emph{absent} edges. Let us
call a $k$-multiset $D$ on $[m]$ \emph{bad} if $D\not\in E$
and $D\not=\multiset{\rep{i}{k}}$ for some $i\in R$. Let 
 $$
 B:=\big(G\setminus \blow{E}{V_1,\dots,V_m}\big)\setminus 
\left(\cup_{i\in R} 
{V_i\choose k}\right)
 $$
 and call edges in $B$ \emph{bad}. Equivalently, an edge of $G$ is bad if
its profile is bad. Define $a:=|A|$ and
$b:=|B|$. Our aim is to achieve that $a=b=0$.

Our next modification is needed to ensure later that~\req{Bxi} holds. Roughly speaking, we want
a property that the number of bad edges cannot be decreased much if 
we move one vertex between parts. 
Unfortunately, we cannot just take a partition structure $\B V$ that minimises $b$ because then we
do not know how to guarantee the high minimum degree of $H$ (another property important in our proof). 
Nonetheless, we can simultaneously satisfy both properties with some extra work.

Namely, we modify $H$ 
as follows (updating $A$, $B$, $\B V$, etc, as we proceed). 
If there is a vertex $x\in [n]$ such that by moving it to
another part $V_i$ we decrease $b$ by at least 
$c_2{n-1\choose k-1}$, then we pick $y\in V_i$ of maximum $H$-degree
and make $x$ a clone of $y$. (Note that the new value of $b$
depends only on the index $i$ of the part $V_i$ but not on the 
choice of $y\in V_i$.)
Clearly, we perform this operation at most 
$c_1{n\choose k}/c_2{n-1\choose k-1}=c_1n/(c_2k)$ times because we initially 
had $b\le |G\bigtriangleup H|\le c_1{n\choose k}$. Thus, we have
at all steps
of this process (which affects at most $c_1n/(c_2k)$ vertices of $H$) that,
trivially,
 \begin{eqnarray}
  |V_{\B i}| &\ge& 2c_4n -\frac{c_1n}{c_2k}\ \ge\ c_4n,\qquad
\mbox{for all legal
$\B i$ with $|\B i|\le \ell$},\label{eq:PartSizes}\\
 |G\bigtriangleup H|&\le& c_1{n\choose k} + \frac{c_1n}{c_2k} {n-1\choose k-1} 
\ \le\ c_2 {n\choose k}.\label{eq:M0End}
 \end{eqnarray}
 It follows that at every step each part $V_i$ had a vertex
of degree at least $(\Lambda_P-c_2/2){n-1\choose k-1}$ 
for otherwise the edit distance between 
$H$ and $G$ at that moment would be at least $\frac{c_2}3 {n-1\choose k-1}\times
c_4n$ by \req{DeltaG} and \req{PartSizes}, 
contradicting the first inequality in \req{M0End}. This implies
that every time we clone a vertex it has a high degree. Thus  we have 
by \req{MinDegStage1} that, additionally to \req{PartSizes} and~\req{M0End},
the following holds at the end of this process:
 \begin{equation}
 \delta(H) \ge \big(\density{P}-\max(3c_1,c_2/2)\big)
{n-1\choose k-1} - \frac{c_1n}{c_2k} {n-2\choose k-2}\ \ge\ 
(\density{P}-c_2) {n-1\choose k-1}.\label{eq:DeltaH1}
 \end{equation}

If we take the
union of $\blow{E}{V_1,\dots,V_m}$ with $\cup_{i\in R} G[V_i]$, then the
obtained $k$-graph is still $\C F_{M_0}$-free by Lemma~\ref{lm:FmFree} and has 
exactly $a-b+|G|$ edges. The
maximality of $G$ implies that
 \beq{b>a}
 b\ge a.
 \eeq

Suppose that $b>0$ for otherwise $a=b=0$ and the lemma is proved. 
Let
 $$
 H':=H\setminus\left(\cup_{\B i\in R^{\ell}} H[V_{\B i}]\right)
 $$ 
 be obtained from $H$ by ``truncating'' it down to the first $\ell$ levels.

Let us show that the maximal degree of $B$ is small, namely that
 \beq{DeltaB1}
 \Delta(B)< c_3 {n-1\choose k-1}.
 \eeq
 Suppose on the contrary that $d_B(x)\ge c_3 {n-1\choose k-1}$ for 
some $x\in [n]$.

It may be helpful to informally illustrate our argument leading to a
contradiction on
the special case when $P$ is as in \req{example}. 
Assume that Lemma~\ref{lm:tight} returns $\ell=1$ and $F=\blow{E}{W_1,W_2}$ 
with $|W_i|=6$ (although some smaller $W_i$'s will also suffice). 
Here we have $H'=\blow{E}{V_1,V_2}$.
Suppose, for example, that the vertex $x$ contradicting \req{DeltaB1}
is in $V_2$. Let $B_{x,2}:=B_x=G_x\cap\big({V_2\choose 2}\cup {V_1\choose
2}\big)$
be the link of $x$ in the bad 3-graph $B$. Next, 
let $B_{x,1}:=G_x\cap \blow{K_2^2}{V_1,V_2}$ be the set of pairs
that would form a bad edge with $x$ if $x$ is moved to $V_1$. By
our assumption, we have $|B_{x,2}|\ge c_3{n-1\choose 2}$. Also,
$|B_{x,1}|\ge (c_3-c_2){n-1\choose 2}$ for otherwise we would have
moved $x$ to $V_1$, thus decreasing $b$ substantially. Take any $\B D=(D_1,D_2)$,
where $D_i\in B_{x,i}$. Consider arbitrary 6-vertex sets $Z_1\subseteq V_1$
and $Z_2\subseteq V_2\setminus\{x\}$ such that $Z_1\cup Z_2$ contains the set $D_1\cup D_2$. 
It is impossible that $\blow{E}{Z_1,Z_2}\subseteq G$ because the $3$-graph $F_{\B D}$
obtained from $\blow{E}{Z_1,Z_2}$ by adding the edges $D_1\cup\{x\}$ and
$D_2\cup\{x\}$ 
belongs to $\C F_{13}$.
(Indeed, if we could embed $F_{\B D}$ into a $P$-construction, then 
all vertices of $Z_1$ and $Z_2$ would have to go into ``correct'' parts by
the rigidity of $F\cong \blow{E}{Z_1,Z_2}$, leaving no way to
fit $x$.)
Thus $Z_1\cup Z_2$ contains at least one absent edge. Finally, our
lower estimates 
on $|B_{x,i}|$ translate with some work into  a lower bound
on $a$ that contradicts \req{M0End}.

Let us give the general argument. For $i\in[m]$ let the $(k-1)$-graph 
$B_{x,i}$ consist of those $D\in G_x$ such that
if we add $i$ to the profile of $D$ then the obtained $k$-multiset
is bad. In other words, if we move $x$ to $V_i$, then $B_{x,i}$ 
will be the link of $x$ with respect to the (updated) bad $k$-graph $B$.
By the definition of $H$, we have  
 \beq{Bxi}
 |B_{x,i}|\ge 
(c_3-c_2){n-1\choose k-1},\qquad \mbox{for every $i\in[m]$.}
 \eeq

For $\B D=(D_1,\dots,D_m)\in \prod_{i=1}^m B_{x,i}$
let $F_{\B D}$ be the $k$-graph that is constructed as follows. Recall that $F$ is the rigid
$P$-construction
given by Lemma~\ref{lm:tight} and $\B W$ is its
partition structure.
By relabelling vertices of $F$, we can assume that $x\not\in V(F)$ while 
$D:=\cup_{i=1}^m D_i$ is a subset 
of $V(F)$ so that for every $y\in D$ we have
$\branch{F}{y}=\branch{H'}{y}$, that is, $y$ has the
same branches in both $F$ and $H'$. (Note 
that both $k$-graphs have the same maximum height $\ell$.) This is possible
because
each part of $F$ of height at most $\ell$ has at least $m(k-1)\ge |D|$
vertices. Finally, add $x$ as a new vertex 
and the sets $D_i\cup\{x\}$
for $i\in [m]$ as edges, obtaining the $k$-graph $F_{\B D}$. 

\bcl{FBD}{For every $\B D\in 
\prod_{i=1}^m B_{x,i}$ we have $F_{\B D}\in \C F_{\infty}$.}

\bcpf Suppose on the contrary that we have an embedding $f$ of $F_{\B D}$
into some $P$-construction with the partition structure $\B U$. By
the rigidity of $F$, we can assume that
$f(W_i)\subseteq U_i$ for every $i\in[m]$. Let $i\in[m]$ satisfy $f(x)\in U_i$.
But then the edge $D_i\cup\{x\}\in F_{\B D}$ is mapped
into a non-edge because $f(D_i\cup\{x\})$ has bad profile
with respect $U_1,\dots,U_m$ by the
choice of $D_i\in B_{x,i}$, a contradiction.\ecpf

For every vector $\B D\in \prod_{i=1}^m B_{x,i}$ and every
map $f:V(F_{\B D})\to V(G)$ such that $f$ is the identity on
$D\cup\{x\}$ and $f$ preserves branches of height up to $\ell$ on all
other vertices, the image $f(F_{\B D})$ has to contain some $X\in
\OO G$ by Claim~\ref{cl:FBD}. (Note that $G$ is $F_{\B D}$-free
since $v(F_{\B D})\le M_0$.) Also,
 $$
 f\big(F_{\B D}\setminus\{D_1\cup\{x\},\dots,D_m\cup\{x\}\}\big)\subseteq H',
 $$
  that is, the ``base'' copy  of $F$ on which $F^{\B D}$ was built is embedded by $f$ into $H'$. On the other hand,
 each of the edges $D_1\cup\{x\},\dots,D_m\cup\{x\}$ of 
$F_{\B D}$ that contain $x$
is mapped to an edge of $G$ (to itself). Thus $X\in 
H'\setminus G$ and $X\not\ni x$. Any such $X$ can
appear, very roughly, for at most ${w\choose k-1}^m\, (w+1)!\, n^{w-k}$ choices
of $(\B D,f)$, where 
$w:=v(F)=v(F_{\B D})-1$.  On the other hand, the total number of choices of
$(\B D, f)$ is at least $\prod_{i=1}^m |B_{x,i}|\ge 
\big((c_3-c_2){n-1\choose k-1}\big)^m$  times
$(c_4n/2)^{w-(k-1)m}$ (since every part of
$H'$ has at least $c_4n$ vertices by~\req{PartSizes}). We conclude that
 $$
 |H\setminus G|\ge |H'\setminus G|\ge 
\frac{\left((c_3-c_2){n-1\choose k-1}\right)^m
\times(c_4n/2)^{w-(k-1)m}}{{w\choose k-1}^m\,(w+1)!\, n^{w-k}}>c_2{n\choose k}.
 $$
 However, this contradicts~\req{M0End}. Thus~\req{DeltaB1} is proved.

Next, we show (in Claim~\ref{cl:M1Deg} below) that \emph{every} bad edge
$D$ 
intersects
$\Omega(c_3n^{k-1})$
absent edges. Again, let us first illustrate the
proof on the case when $P$ is as in \req{example}.
Suppose that $D=\{y_1,y_2,z\}$ with, say, $y_1,y_2\in V_1$ and $z\in V_2$. For $i=1,2$, 
let $D_{y_i}$ be an edge of $G[V_1]$ such that $D_{y_i}\cap D=\{y_i\}$; there
are many such edges because we have by Part~\ref{it:OmegaNDelta} of
Lemma~\ref{lm:OmegaN} that, for example,
 \beq{cl:Level1Deg}
  d_{G[V_i]}(y)\ge c_4{n-1\choose k-1},\quad \mbox{for all $i\in R$ and $y\in
V_i$}.
 \eeq
  Let $\B
D=(D_{y_1},D_{y_2})$.
Take arbitrary $6$-sets $Z_i\subseteq V_i\setminus D$ such that
$Z_1\supseteq (D_{y_1}\cup D_{y_2})\setminus D$. Let $F^{\B D}$ be the
$3$-graph obtained
from $\blow{E}{Z_1,Z_2\cup\{z\}}$ by adding vertices $y_1,y_2$ and
edges
$D,D_{y_1},D_{y_2}$.
It is not hard to show that $F^{\B D}$ belongs to $\C F_{15}$. Thus
$\blow{E}{Z_1,Z_2\cup\{z\}}\not\subseteq G$ and we arrive at some $Y\in A$.
It is impossible that the obtained absent edge $Y$ is
disjoint from $D$ for at least half of choices of $(D_{y_1},D_{y_2},Z_1,Z_2)$,
for otherwise $|A|$ is too large. If $Y\cap D$ is not empty, then it
consists of the unique vertex $z$. Some counting gives 
the desired lower bound on the number of absent edges intersecting $D$. 
Note that this counting would not work if the obtained $Y\in A$ could
share more than one vertex with $D$. This is the reason why 
we do not allow $Z_1\cup Z_2$ to share more than one vertex with $D$;
we make these sets disjoint in the general proof  for the notational
convenience.

Let us present the general argument. Let $D\in B$ be an arbitrary bad edge. 
For each $i\in R$ and $y\in D\cap V_i$ pick
some $D_y\in G[V_i]$ such that $D_y\cap D=\{y\}$; it exists
by~\req{cl:Level1Deg}. Let 
$\B D:=(D,\{D_y\mid y\in D\cap V_R\})$. (Recall that
$V_R=\cup_{i\in R} V_i$.) 
We define the $k$-graph
$F^{\B D}$ 
using the rigid $k$-graph $F$ as follows. By re-labelling $V(F)$, we can assume
that $X\subseteq V(F)$,
where 
 \beq{X}
 X:=\cup_{y\in D\cap V_R} D_y\setminus\{y\},
 \eeq
so that for every $x\in X$ its branches in $F$ and $H'$
coincide. Again, there is enough space inside $F$ to accommodate all
$|X|\le k(k-1)$ vertices of $X$. Assume also that $D$ is disjoint 
from $V(F)$. The vertex
set of $F^{\B D}$ is $V(F)\cup D$. Starting with the edge-set of $F$, add
$D$ and each $D_y$ with $y\in D\cap V_R$.
Finally, for every $y\in D\cap V_i$ with $i\in [m]\setminus R$ pick
some $z\in W_i$ and add $\{Z\cup \{y\}\mid Z\in F_z\}$ to the edge set,
obtaining the $k$-graph $F^{\B D}$.
The last step can be viewed as enlarging the part $W_i$ by $D\cap V_i$
and adding those edges that are stipulated by the pattern $P$ and
intersect $D$ in at most one vertex.

\bcl{FBD2}{For every $\B D$ as above, we have $F^{\B D}\in\C F_\infty$.}

\bcpf Suppose on the contrary that we have an embedding $f$ of $F^{\B D}$
into some $P$-construction with the partition structure $\B U$. We
can assume by the rigidity of $F$, that $f(W_i)\subseteq U_i$
for each $i$. 

Take
$y\in D\cap V_i$ with $i\in R$. 
The $(k-1)$-set $f(D_y\setminus\{y\})$ lies entirely
inside $U_i$. We cannot have $f(y)\in U_j$ with $j\not=i$ because
otherwise the profile of the edge $f(D_y)$ is $\multiset{\rep{i}{k-1},j}$,
contradicting Lemma~\ref{lm:density=1}.  Thus $f(y)\in U_i$.

Next, take any $y\in D\cap V_i$
with $i\in [m]\setminus R$. Pick some $z\in W_i$.
By the rigidity of $F$, if
we fix the restriction of $f$ to 
$V(F)\setminus\{z\}$, then $U_i$ is the only part where $z$ can be mapped to. 
By definition,
$y$ and $z$ have the same link $(k-1)$-graphs in  $F^{\B D}$ when restricted to
$V(F)\setminus\{y,z\}$. Hence, $f(y)$ necessarily belongs to $W_i$.

Thus the edge $f(D)$ has the same profile as $D\in B$,
a contradiction.\ecpf

\bcl{M1Deg}{For every $D\in B$ there are at least $kc_3{n-1\choose k-1}$
absent edges $Y\in A$ with $|D\cap Y|=1$.}

\bcpf Given $D$ choose the sets $D_y$, for $y\in D\cap V_R$, as before 
Claim~\ref{cl:FBD2}. The condition $D_y\cap D=\{y\}$ rules out
at most $k {n-2\choose k-2}$ edges for this $y$. Thus
by~\req{cl:Level1Deg} there are, for example, at least
$(c_4/2){n-1\choose k-1}$ choices of each $D_y$. Form the $k$-graph
$F^{\B D}$ as above and consider potential injective embeddings $f$ of 
$F^{\B D}$ into
$G$ that are the identity on $D\cup X$ and map every other 
vertex of $F$ into a vertex of $H'$ with the same
branch, where $X$ is defined by \req{X}. For every vertex $x\not\in D\cup X$ we have at least $c_4n/2$
choices for $f(x)$ by \req{PartSizes}. By Claim~\ref{cl:FBD2}, $G$ does not contain $F^{\B D}$
as a subgraph so its image under $f$ contains some $Y\in\OO G$.
Since $f$ maps $D$ and each $D_y$ to an edge of $G$
(to itself) and 
 $$
 f\left(F^{\B D}\setminus(\{D\}\cup \{D_y\mid y\in D\cap V_R\})\right)\subseteq H',
 $$ 
 we have that $Y\in H'$. 
The number of choices of $(\B D,f)$ is at least
 $$
 \left((c_4/2){n-1\choose k-1}\right)^{|D\cap V_R|}\times
(c_4n/2)^{w-(k-1)|D\cap V_R|}
\ge (c_4n/4k)^{w},
 $$
 where $w:=v(F)$. Assume that for at least half of the time the obtained set
$Y$ intersects $D$ for otherwise we get a contradiction to~\req{M0End}:
 $$
 |H'\setminus G|\ge \frac12 \times \frac{(c_4n/4k)^{w}}{{w\choose k-1}^k\,
(w+k)!\, n^{w-k}}>c_2{n\choose k}.
 $$

By the definitions of 
$F^{\B D}$ and $f$, we have that $|Y\cap D|=1$ and
$Y\in A$. 
Each such $Y\in A$ is counted for at most ${w\choose k-1}^k\,(w+k)!\,
n^{w-k+1}$ 
choices of $f$. Thus the number of such $Y$ is at least
$\frac1{2}(c_4n/4k)^{w}/({w\choose k-1}^k\,(w+k)!\, n^{w-k+1})$, implying the
claim.\ecpf

Let us count the number of pairs $(Y,D)$ where $Y\in A$, $D\in B$,
and $|Y\cap D|=1$. On one hand, each bad edge $D\in B$ creates at least
$kc_3{n-1\choose k-1}$ such pairs by Claim~\ref{cl:M1Deg}.
On the other hand, we trivially
have at most $a k\Delta(B)$ such pairs. By~\req{b>a}, we have $b kc_3{n-1\choose k-1}\le a k\Delta(B)\le b k\Delta(B)$. 
Since $b\not=0$, we obtain a contradiction to~\req{DeltaB1}.
This finishes the proof of Lemma~\ref{lm:Max}.\epf

Let us state a special case of a result of R\"odl and
Schacht~\cite[Theorem~6]{rodl+schacht:09} that we will need.
  
\begin{lemma}[Strong Removal Lemma~\cite{rodl+schacht:09}]\label{lm:RS} For every $k$-graph family $\C F$ 
and $\e>0$ there are $\delta>0$, $m$, and $n_0$ such 
that the following holds. Let $G$ be a $k$-graph on $n\ge n_0$ vertices 
such that for every $F\in \C F$ with $v(F)\le m$ the number of $F$-subgraphs in $G$
is at most $\delta n^{v(F)}$. Then $G$ can be made $\C F$-free
by removing at most $\e {n\choose k}$ edges.\qed
\end{lemma}

\begin{lemma}\label{lm:edit} For every $c_0>0$ there is $M_1$ such that
every maximum $\C F_{M_1}$-free $G$ with $n\ge M_1$ vertices is 
$c_0{n\choose k}$-close
to a $P$-construction.\end{lemma}
 \bpf Lemma~\ref{lm:RS} 
gives $M_1$ such that any $\C F_{M_1}$-free $k$-graph $G$
on $n\ge M_1$ vertices can be made into an $\C F_\infty$-free 
$k$-graph $G'$ by removing at most $c_0{n\choose k}/2$ edges. By 
Lemma~\ref{lm:FInfty}, $G'$ embeds into some $P$-construction $H$
with $v(H)=v(G')$.
Assume that $V(H)=V(G')$ and the identity map is an embedding of 
$G'$ into~$H$.

Since $H$ is $\C F_{M_1}$-free, 
the maximality of $G$ implies that $|G|\ge |H|$. Thus 
$|H\setminus G'|\le c_0{n\choose k}/2$ and we can transform $G'$ into
$H$ by changing at most $c_0{n\choose k}/2$ further edges.\epf

\subsection{Proof of Theorem~\ref{th:Fin}: Putting All Together}

We are ready to prove Theorem~\ref{th:Fin}. It is trivially true if $\density{P}=0$ or $1$ by the discussion
after Lemma~\ref{lm:density=1}, so we can
assume that $P$ is proper. 
Apply Lemma~\ref{lm:Max} which
returns $c_0$ and $M_0$. Next, Lemma~\ref{lm:edit} on input $c_0$ returns
some $M_1$. 

Let us show that $M=\max(M_0,M_1)$ works in
Theorem~\ref{th:Fin}. We use induction on $n$.  Let $G$ be any maximum
$\C F_M$-free $k$-graph on $[n]$. Suppose that $n> M$ for otherwise we
are done by Lemma~\ref{lm:FInfty}. Thus Lemma~\ref{lm:edit} applies
and shows that $G$ is $c_0{n\choose k}$-close to some
$P$-construction. Lemma~\ref{lm:Max} returns a partition
$[n]=V_1\cup\dots\cup V_m$ such that \req{BottomLevel} holds. 

Let $i\in R$ be arbitrary. By Lemma~\ref{lm:FmFree} if we replace 
$G[V_i]$ by a maximum $\C F_M$-free $k$-graph, then the new
$k$-graph on $V$ is still $\C F_M$-free.
By the maximality of $G$, we conclude
that $G[V_i]$ is a maximum $\C F_M$-free $k$-graph. By the induction
hypothesis (note that $|V_i|\le n-1$), $G[V_i]$ is a
$P$-construction. 

It follows that $G$ is a $P$-construction itself, which implies all
claims of Theorem~\ref{th:Fin}.\qed

\section{Proof of Theorem~\ref{th:irrat}}\label{irrat}

By Corollary~\ref{cr:Fin} it is enough to exhibit, for every
$k\ge 3$, a pattern $P$ such that $\density{P}$ is irrational.

Given $k\ge 3$, 
let $\ell$ be any prime number that does not divide $k$ such that 
$2\le \ell< k$. If $k$ is odd, we can take $\ell=2$. For 
even $k$ we can take $\ell$ to be any prime with $k/2<\ell< k$;
it exists by Bertrand's postulate. Take
$P=(2,E,\{1\})$, where $E$ consists of the single multiset $\multiset{\rep{1}{k-\ell},\rep{2}{\ell}}$. In other words, a $P$-construction
on $V$ is obtained by partitioning $V=V_1\cup V_2$ with $V_1\not=V$, adding
all $k$-sets that intersect $V_1$ in exactly $k-\ell$ vertices, and doing
recursion inside $V_1$. If $k=3$, we get the familiar pattern from
\req{example}.

Let $\B x=(x_1,x_2)\in\C X$ be an optimal vector. We have 
by~\req{opt} that $\density{P}/{k\choose \ell}=r(x_1)$, where
 $$
 r(x):=\frac{(1-x)^\ell x^{k-\ell}}{1-x^k}.
 $$
 By Part~\ref{it:f} of Lemma~\ref{lm:Omega1} 
the real $x_1$ maximises $r$ over $(0,1)$. Thus $x_1$ is a root of the
derivative of $r$. We
have 
 $$
 r'(x)=-\frac{(1-x)^{\ell}x^{k-\ell-1}}{(1-x^k)^2}\, g(x),
 $$ 
 where $ g(x):=
\ell (x^{k-1}+\dots+1)-k$.
 Since $x_1\not=0,1$, it is a root of $g$. 

The polynomial $g(x)$ is
irreducible (by Eisenstein's criterion). Indeed, if we can factorise $g(x)=(a_mx^m+\dots+a_0)(b_{k-1-m}x^{k-1-m}+\dots+b_0)$ in $\I Z[x]$, 
then exactly
one of $a_m$ and $b_{k-1-m}$ is divisible by $\ell$, say $b_{k-1-m}$ is.
Since $b_0 a_0=\ell-k$, we have that $\ell$ does not divide $b_0$. Thus there
is $i$ such that $\ell$ does not divide $b_i$ but $\ell$ divides each $b_j$ with $j>i$. 
But then the coefficient at $x^{m+i}$ is congruent to $a_mb_i$ modulo $\ell$,
which implies that $m=0$ and the first factor is equal to $\pm1$.

Thus $x_1$ is irrational but we still have to show that $\density{P}$ is
irrational. Suppose on the contrary that $\density{P}=s/t$ with $s,t\in\I Z$.
Note that $(x-1)g(x)=\ell x^k-kx+k-\ell$. Thus
$x_1^k=(kx_1-k+\ell)/\ell$.  Substituting this in $\Lambda_P={k\choose
\ell}r(x_1)$, we infer
that 
 $$
 \frac st={k\choose \ell}\frac{(1-x_1)^\ell
x_1^{k-\ell}}{1-(kx_1-k+\ell)/\ell}
= {k-1\choose \ell-1}(1-x_1)^{\ell-1}x_1^{k-\ell}.
 $$
 Thus $x_1$ is a root of the polynomial 
 $$
 h(x):=t\,(k-1)!\,(1-x)^{\ell-1} 
x^{k-\ell}-s\,(\ell-1)!\,(k-\ell)!\ \in\ \I Z[x]
 $$ 
 which has to be divisible by the irreducible polynomial $g$. Since these
polynomials have the same degree $k-1$, $h$ is a constant multiple of
$g$. But the highest two coefficients of $g$ are the same while those of $h$
have different signs, a contradiction. Thus $\density{P}$ is irrational,
proving Theorem~\ref{th:irrat}.\qed

\section{Proof of Proposition~\ref{pr:closed}}\label{closed}

Here we prove Proposition~\ref{pr:closed}. The proof is motivated
by the
emerging theory of the limits of discrete structures, 
see e.g.~\cite{BCLSV:08,elek+szegedy:12,lovasz:lngl,lovasz+szegedy:06}.
Also, an intermediate result that we obtain (Theorem~\ref{th:CT}) 
may be of independent interest. 
We make the presentation essentially self-contained by
restricting ourselves to only one aspect of hypergraph limits.
In particular, we do not rely on the machinery developed by Elek and 
Szegedy~\cite{elek+szegedy:12}.

Let $F$ and $G$ be $k$-graphs. A \emph{homomorphism} from $F$ to $G$ is a map
$f:V(F)\to V(G)$, not necessarily injective, such that $f(A)\in G$
for every $A\in F$. (Thus, embeddings are precisely injective
homomorphisms.) Let the \emph{homomorphism density} $t(F,G)$ be the probability
that a random map $V(F)\to V(G)$, with all $v(G)^{v(F)}$ 
choices being equally likely, is a homomorphism. For
example, we have $t(K_k^k,G)=k!\,|G|/v(G)^k$.

Let $\CG^{(k)}$
consist of all $k$-graphs up to isomorphism.
A sequence $(G_i)_{i=1}^\infty$
of $k$-graphs \emph{converges to a function} $\phi:\CG^{(k)}\to [0,1]$ if 
the sequence is
\emph{increasing} (i.e.\ $v(G_1)<v(G_2)<\dots$)
and for every $k$-graph $F$ we have
$\lim_{i\to\infty} t(F,G_i)=\phi(F)$. Clearly,
the convergence is not affected if we modify
$o(v(G_i)^k)$ edges in each $G_i$. Let $\LIM{k}$ consist of all 
possible functions $\phi$ that can be obtained in the above manner.

Given a family $\C F$ of forbidden $k$-graphs, let 
$\C T(\C F)\subseteq \LIM{k}$ consist of all possible limits of increasing
sequences $(G_i)_{i=1}^\infty$ such that 
$\lim_{i\to\infty} \rho(G_i)=\pi(\C F)$
and each $G_i$ is $\C F$-free. In other words,
$\C T(\C F)$ is the set of the limits of almost maximum
$\C F$-free $k$-graphs. The standard
diagonalisation argument shows that every increasing
sequence has a convergent subsequence; in particular,
$\C T(\C F)\not=\emptyset$. Let $\C T^{(k)}$ be the union
of $\C T(\C F)$ over all $k$-graph families $\C F$. We have
 \beq{Pi=ImCT}
 \InfPi{k}=\{\phi(K_k^k)\mid \phi\in \C T^{(k)}\}.
 \eeq

Let the \emph{blow-up closure $\cl{\C F}$} of $\C F\subseteq \CG^{(k)}$
consist of all $k$-graphs $F$ such that some blow-up of $F$
is not $\C F$-free. Clearly, $\C F\subseteq \cl{\C F}$. Also,
it is easy to see that by applying the blow-up closure twice we
get the same family $\cl{\C F}$. 

\begin{lemma}\label{lm:T(CF)} For every $\C F\subseteq \CG^{(k)}$ 
and $\e>0$ there is $n_0$ such that any $\C F$-free $k$-graph
$G$ with $n\ge n_0$ vertices can be made $\cl{\C F}$-free by removing
at most $\e {n\choose k}$ edges.

In particular, it follows that $\pi(\C F)=\pi(\cl{\C F})$ and $
\C T(\C F)=\C T(\OO{\C F})$. 
\end{lemma}
 \bpf Let Lemma~\ref{lm:RS} on input $(\C F,\e)$ return
$m$ and $\delta>0$. Let $n$ be large and $G$
be an arbitrary $\C F$-free $k$-graph on $[n]$.
For each $F\in \cl{\C F}$ there is $s$ such that
$G$ is $F(s)$-free. As it is well known (see 
e.g.~\cite[Theorem~3]{brown+simonovits:84}), $G$
contains at most $\delta n^{v(F)}$ copies of $F$ 
for all large $n$. 
Since there are only finitely many non-isomorphic $k$-graphs
on at most $m$ vertices, we can satisfy the above bound 
for all such $k$-graphs by taking $n$ large. Now Lemma~\ref{lm:RS} applies,
giving the first part.

Thus any increasing sequence $(G_i)_{i=1}^\infty$ of asymptotically maximum $\C
F$-free $k$-graphs can be converted into that for $\OO{\C F}$ by modifying
$o(v(G_i)^k)$ edges in each $G_i$. This modification does not affect,
for any fixed $F$,
the limit of $t(F,G_i)$ as $i\to\infty$, implying the second part.\epf

Recall
that the Lagrangian of a $k$-graph $G$ on $[n]$ is 
$\Lambda_G=\max\{\lambda_G(\B x)\mid \B x\in
\IS_n\}$; equivalently, $\Lambda_G$ is the Lagrangian $\Lambda_P$ of the pattern 
$P:=(n,G,\emptyset)$ as defined by~\req{LagrangianP}. 
We have the following characterisation of the set $\C T^{(k)}$.

\begin{theorem}\label{th:CT} For $\phi\in\LIM{k}$,
the following are equivalent:
 \begin{enumerate}
  \item\label{it:CT1} $\phi\in \C T^{(k)}$;
 \item\label{it:CT2} $\phi$ is a limit of an increasing $k$-graph sequence
$(G_i)_{i=1}^\infty$
such that $\rho(G_i)-\Lambda_{G_i}\to 0$;
  \item\label{it:CT3} $\phi(F)=0$ for every $k$-graph $F$ with
$t(K_k^k,F)>\phi(K_k^k)$.
 \end{enumerate}
\end{theorem}
\bpf \textbf{1) $\Rightarrow$ 2):} Let $\phi\in\C T^{(k)}$, say
$\phi\in\C T(\C F)$. By Lemma~\ref{lm:T(CF)} we can assume
that $\cl{\C F}=\C F$. Let $(G_i)_{i=1}^\infty$ be
a sequence of almost maximum $\C F$-free $k$-graphs
that converges to $\phi$. Take any $i$ and let $n:=v(G_i)$. 
Since $\cl{\C F}=\C F$, any blow-up
of $G_i$ is still $\C F$-free. Also, the limit superior of 
the edge densities
attained by increasing blow-ups of $G_i$ is exactly $\Lambda_{G_i}$. Thus
$\Lambda_{G_i}\le \pi(\C F)=\rho(G_i)+o(1)$.
On the other hand, we have $\Lambda_{G_i}\ge \lambda_{G_i}(1/n,\dots,1/n)= k!\, |G_i| /n^k$, giving the converse inequality.\medskip

\noindent\textbf{2) $\Rightarrow$ 3):} Let $\phi$ satisfy~2). Suppose on the contrary that some $F$ on $[n]$ violates~3). 
Pick a sequence $(G_i)_{i=1}^\infty$ given by~2). By the definition of convergence,
$G_i$ contains $F$ as a subgraph for all large $i$. But then
 $$
 \Lambda_{G_i}\ge \Lambda_{F}\ge \frac{k!\, |F|}{n^k} =t(K_k^k,F)\ge
\phi(K_k^k)+\Omega(1),
 $$ 
 contradicting $\phi(K_k^k)=\lim_{i\to\infty}\rho(G_i)=\lim_{i\to\infty}\Lambda_{G_i}$.\medskip

\noindent\textbf{3) $\Rightarrow$ 1):} Given $\phi$ as in~3), 
let 
 $$
 \C F:=\{F\in\CG^{(k)}\mid \phi(F)=0\}.
 $$

Let $H_n$ be a maximum $\C F$-free $k$-graph on $[n]$. Since
$H_n\not\in \C F$, we have $\phi(H_n)>0$. Thus, by 3),
 $$
\phi(K_k^k)\ge t(K_k^k,H_n)=\frac{k!\,|H_n|}{n^k}\ge\rho(H_n)+O(1/n).
 $$
 By letting $n\to\infty$, we obtain that
 \beq{pi=phi}
 \phi(K_k^k)\ge \pi(\C F).
 \end{equation}

Let us show that $\phi\in \C T(\C F)$. Take any increasing sequence
$(G_i)_{i=1}^\infty$ convergent to $\phi$. Let $F\in \C F$.
Since $\phi(F)=0$, the number of
$F$-subgraphs 
in $G_i$ is at most
$o(v(G_i)^{v(F)})$.
By Lemma~\ref{lm:RS},
we can remove $o(v(G_i)^k)$
edges in each $G_i$, obtaining an $\C F$-free $k$-graph $G_i'$. Thus 
 $$
 \pi(\C F)\ge \lim_{i\to\infty} \rho(G_i')=\phi(K_k^k)
 $$
 and, by \req{pi=phi}, this is equality. The obtained sequence
$(G_i')_{i=1}^\infty$
of almost maximum 
$\C F$-free $k$-graphs
still converges to $\phi$. This shows that $\phi\in\C T(\C F)\subseteq\C
T^{(k)}$.\epf

Let us view $\LIM{k}$
as a subset of $[0,1]^{\CG^{(k)}}$, where the latter set is endowed with
the product (or pointwise convergence) topology. If we
identify each $k$-graph $G$ with the sequence $(t(F,G))_{F\in\CG^{(k)}}$,
then this topology gives exactly the above convergence. Moreover,
the set $\LIM{k}$, as the topological closure of $\CG^{(k)}$,
is a closed subset of $[0,1]^{\CG^{(k)}}$.

\begin{corollary}\label{cr:CT} For every $k\ge 2$ the set $\C T^{(k)}$ is 
a closed subset of $[0,1]^{\CG^{(k)}}$.
\end{corollary}
 \bpf The third characterisation
of Theorem~\ref{th:CT} shows that
 \beq{CT}
 \C T^{(k)}=\cap_{F\in\CG^{(k)}} \left( \big\{\phi\in \LIM{k}\mid \phi(F)=0\big\}
\cup \big\{\phi\in \LIM{k}\mid \phi(F)\ge t(K_k^k,F)\big\} \right).
 \eeq
 For every
$F\in\CG^{(k)}$, 
the map $\phi\mapsto \phi(F)$ is continuous 
(it is just the projection of $[0,1]^{\CG^{(k)}}$ onto
the $F$-th coordinate). We see by~\req{CT} that $\C T^{(k)}$, 
as the intersection
of closed sets, is itself closed.\epf

\bpf[Proof of Proposition~\ref{pr:closed}.]
Let $a_i\in \InfPi{k}$ with $a_i\to a$ as $i\to\infty$.
By Part~\ref{it:CT2} of Theorem~\ref{th:CT} 
we can find, for each $i\in\I N$, a $k$-graph $H_i$ such that $v(H_i)>i$ and both
$\Lambda_{H_i}$ and edge density of $H_i$ are within $1/i$ from $a_i$. 
By passing to a subsequence
we can additionally assume that the $k$-graphs 
$H_i$ converge to some $\phi\in\LIM{k}$.
This $\phi$ satisfies Part~2 of 
Theorem~\ref{th:CT} and thus belongs to $\C T^{(k)}$. Thus
$a=\phi(K_k^k)$ belongs to $\InfPi{k}$, as required.

Alternatively, by Tychonoff's theorem, $[0,1]^{\CG^{(k)}}$
is compact. By Corollary~\ref{cr:CT}, $\C T^{(k)}$
is compact. By~\req{Pi=ImCT}, 
$\InfPi{k}$ is a continuous image 
of $\C T^{(k)}$, so it is compact too. Hence
$\InfPi{k}\subseteq [0,1]$ is closed.\epf

\section{Proof of Theorem~\ref{th:uncountable}}\label{uncountable}

Let $k\ge 3$. Let $\alpha<1$ be a non-jump for $k$-graphs 
(that is, $(\alpha,\alpha+\e)\cap \InfPi{k}\not=\emptyset$ for every
$\e>0$). It exists by the result
of Frankl and R\"odl~\cite{frankl+rodl:84}.  Pick $m$ so that
$\gamma:=\Lambda_{K_m^k}>\alpha$. 
(Such $m$ exists as the assignment $x_i=1/m$
shows that $\Lambda_{K_m^k}\ge k! {m\choose k}/m^k$, which
tends to $1$ as $m\to\infty$.) Let $\tau:=\gamma/(k(m-1))$.
By Part~2 of Theorem~\ref{th:CT}, we can pick, inductively for $i=1,2,\dots$, 
a $k$-graph $H_i$ such that $\beta_i:=\density{H_i}$ belongs to
$(\alpha,\gamma)$ and
 \beq{li}
 0< \beta_{i}-\alpha <
(\beta_{i-1}-\beta_{i})\,\tau^{2k},\qquad
\mbox{for all $i\ge 2$}.
 \eeq 
 Informally speaking, we require that $\beta_1>\beta_2>\dots$ tend to $\alpha$ 
rather fast.

Next, we introduce a new concept that is similar to that of a
$P$-construction. 
Namely, for an infinite set $A=\{a_1<a_2<\dots\}\subseteq \I N$, an
\emph{$A$-configuration} is a $k$-graph $G$ that can be recursively obtained
as follows. Take a partition
$V=V_1\cup\dots\cup V_m$ of the vertex set and add $\blow{K_m^k}{V_1,\dots,V_m}$
to the edge 
set (that is, add all $k$-sets that intersect every part in at most
one vertex).
Inside $V_2$ put some blow-up of $H_{a_1}$. Inside $V_1$ put any
$A'$-configuration, where $A':=A\setminus\{a_1\}$. 

Note that we allow
a part to be everything, e.g.\ we allow $V_1=V$. Let $p_{A,n}$ be the
maximum size of an $A$-configuration on $n$ vertices.
Let $\C F_A$ consist of all $k$-graphs that do not embed into
an $A$-configuration. It is routine to see that Lemmas~\ref{lm:constr}--\ref{lm:regular},
with the obvious modifications, apply to $A$-configurations as well.
In particular, we have $\ex(n,\C F_A)=p_{A,n}$ for all $n$.
Let $\density{A}$ be the limit of $p_{A,n}/{n\choose k}$ as $n\to\infty$;
averaging shows
that this 
ratio is non-increasing (cf Lemma~\ref{lm:lim}). Thus $\density{A}=\pi(\C F_A)$.

In order to show that $|\InfPi{k}|\ge 2^{\aleph_0}$ it is enough to show 
that $\density{A}\not=\density{B}$ 
for every pair of infinite distinct sets $A,B\subseteq \I N$.  We prove the
stronger claim that $\density{A}>\density{B}$ provided
 \beq{MinA-B}
 \min A\setminus B<\min B\setminus A,
 \eeq
 where we agree that $\min X=\infty$ if $X$ is empty.

Let $A=\{a_1<a_2<\dots\}$, 
$B=\{b_1<b_2<\dots\}$, and $\min A\setminus B=a_i$. 
Fix large $\ell$ and let $n\to\infty$. Take a maximum $B$-configuration $G$.
Let $\B V$ be its partition structure, defined in the obvious way. (For
example, 
every index in $\B V$ is of the form $(\rep{1}{j},s)$ for some $j\ge 0$ and
$s\in [m]$.)

If, for some $j\le i$ and infinitely many $n$, the part $V_{\rep{1}{j-1},2}$
(that is, the second part
of the $j$-th level of $G$) is empty, we remove this $b_j$ from $B$.
Clearly, the $k$-graph $G$ remains a maximum $B$-configuration. 
Also, this does not violate
\req{MinA-B}. Thus, by passing to a subsequence of $n$, we can 
assume that $V_{\rep{1}{j-1},2}\not=\emptyset$ for all $j\le i$.
Furthermore, by relabelling parts (if needed), 
we can assume that for every $j\ge 0$
 \beq{12}
 \min\left(|V_{\rep{1}{j},1}|,|V_{\rep{1}{j},2}|\right)
\ge \max\left\{|V_{\rep{1}{j},h}|\mid h=3,\dots,m\right\}.
 \eeq

Let us show by induction on $j=1,\dots,\ell$ that 
 \beq{U12}
 \min\big(|U_1|,|U_2|\big)\ge  \tau\, |U|,
 \eeq
 where $U:= V_{\rep{1}{j-1}}$ and $U_h:=V_{\rep{1}{j-1},h}$ for $h\in[m]$.
Since $G':=G[U]$ is a maximum $\{b_j,b_{j+1},\dots\}$-configuration on
at least $\tau^{j-1}n$ vertices, its edge density is, for example, at
least $\gamma+o(1)$. The argument of Lemma~\ref{lm:regular}
shows that $\delta(G')\ge (\gamma+o(1)) {|U|-1\choose k-1}$. 
It is impossible that $U_2=U$ for otherwise $\Lambda_{G'}\le 
\Lambda_{H_{b_j}}< \gamma$ and $G'$ cannot be maximum
for large $n$. This, our
assumption that $U_2\not=\emptyset$ and~\req{12} imply 
that both $U_1$ and $U_2$ are non-empty. Let $h=1$ or $2$. It is
impossible that $|U_h|>(1-\gamma/k)|U|$
for otherwise a vertex $x$ of  
$U_{3-h}\not=\emptyset$ has too
small $G'$-degree as every edge of $G_x'$ has at least
one other vertex in $U\setminus U_h$. By~\req{12}
we get that $|U_{3-h}|\ge (|U|-|U_{h}|)/(m-1) \ge \tau\,|U|$.
This proves~\req{U12}.

Recall that $i\in\I N$ is defined by $a_i=\min A\setminus B$. Let
$U_1\cup\dots\cup U_m$ be
the partition of $U:=V_{\rep{1}{i-1}}$ in the $B$-configuration $G$.
Let $G'$ be obtained
from $G$ by replacing $G[U_2]$ with a maximum blow-up of $H_{a_i}$
(instead of $H_{b_i}$)
and replacing $G[U_1]$ with the $\{a_j\mid j>i\}$-configuration that has the
same partition structure as the $\{b_j\mid j>i\}$-configuration 
$G[U_1]$. Clearly, $G'$ is an $A$-configuration.
Since $|U_2|\ge \tau|U|$ by \req{U12}, 
the change inside $U_2$ increases the number of edges 
by at least $(\beta_{a_i}-\beta_{a_i+1})\tau^k{u\choose k}+o(n^k)$, where $u:=|U|$.
On the other hand, when we modify $G[U_1]$, 
we replace, for $j>i$, a blow-up of $H_{b_j}$
by another blow-up whose density is at least 
$\alpha+o(1)$. Let $n_j$ be the number of vertices in this part. By \req{U12},
$n_j\le (1-\tau)^{j-i}|U|$ for all $j\le \ell$.
Thus
 \begin{eqnarray*}
 \big|G[U_1]\big|-\big|G'[U_1]\big|&\le &\sum_{j=i+1}^{\ell-1} 
(\beta_{b_j}-\alpha){n_j\choose k}+\sum_{j\ge \ell} {n_j\choose k}+o(n^k)\\
 &\le& (\beta_{a_i+1}-\alpha) {u\choose k} \sum_{j=i+1}^{\ell-1} \left(1-\tau\right)^{(j-i)k}
+{n_{\ell}+n_{\ell+1}+\dots\choose k}+o(n^k)\\
 &\le & \big((\beta_{a_i+1}-\alpha)\tau^{-k} 
+ (1-\tau)^{(\ell-i)k}\big){u\choose k}.
 \end{eqnarray*}
 This is strictly less than $(\beta_{a_i}-\beta_{a_i+1})\tau^k{u\choose k}$
by~\req{li} (and since $\ell=\ell(A,B)$ is large). Thus $|G'|\ge
|G|+\Omega(n^k)$ and indeed $\density{A}>\density{B}$, finishing the proof
of 
Theorem~\ref{th:uncountable}.\qed

\section{Concluding Remarks}\label{conclusion}

If we consider graphs (the case $k=2$), then the Stability Theorem of Erd\H
os~\cite{erdos:67a} and Simonovits~\cite{simonovits:68} answers the question
about the possible asymptotic structure of maximum $\C F$-free graphs.  However,
if we need 
a more precise answer, then the picture is much more complicated
and many questions remain open, including the general \emph{inverse problem} of
describing graphs that are maximum $\C F$-free for
some family $\C F$
(see~e.g.\ \cite{simonovits:68,simonovits:74,simonovits:83:misc}). The
situation with extremal
problems for digraphs and multigraphs is similar (see
e.g.~\cite{brown+erdos+simonovits:85,brown+erdos+simonovits:84,
brown+simonovits:02,sidorenko:93:c}).

Although very few instances of the hypergraph Tur\'an problem have
been solved, there is a variety of constructions giving
best known lower bounds. 
So it is likely that $\FinPi{k}$ contains many further elements
in addition to the values given by Corollary~\ref{cr:Fin}.
 For example, we do not
know if there is a pattern $P$ that gives the same (or better) lower bound $\pi(\{K_5^4\})
\ge \frac{11}{16}$
as the  construction of Giraud~\cite{giraud:90} (see also~\cite{decaen+kreher+wiseman:88}
for generalisations). Roughly speaking, Giraud's construction 
takes an arbitrary 2-colouring of vertices and pairs of vertices 
(with an optimal colouring of pairs being quasi-random)
and decides if a quadruple $X$ is an edge depending on the colouring 
induced by $X$. It would be interesting to decide if
Corollary~\ref{cr:Fin} can be extended
to cover constructions of this type.

In the special case when $E$ consists of simple $k$-sets and $R=\emptyset$,
$\density{P}$ is equal by Lemma~\ref{lm:Omega1} to the well-studied 
Lagrangian of the $k$-graph $E$,
see e.g.~\cite{baber+talbot:12}. Thus
Corollary~\ref{cr:Fin} implies that every value of the Lagrangian belongs to
$\FinPi{k}$, answering a question of 
Baber and Talbot~\cite[Question~29]{baber+talbot:12}. 

One can
show that every proper pattern $P=(m,E,R)$ with
$R\not=\emptyset$ is \emph{complex}, meaning that the number of non-isomorphic
$s$-vertex subgraphs in a large maximum $P$-construction grows faster than
any polynomial of $s$. Indeed, by Lemma~\ref{lm:MaxRigid} for every $\ell$
there is a $P$-construction $F$ with the partition structure $\B V$ 
which is \emph{$\ell$-rigid}, meaning that for every 
$\B i\in R^s$ with $s\le \ell$ the induced $P$-construction
$F[V_{\B i}]$ is rigid. 
 Additionally,
we can assume that $|V_{\B i}|\ge k$ for each legal $\B i$ of length at most $\ell+1$. Thus
if we add any $n-v(F)$ vertices,  the new $k$-graph $F'$ is still $\ell$-rigid
by Lemma~\ref{lm:AddRigid}. There are at least $\ell$ different parts at
the bottom $\ell$ levels for placing these extra vertices. The rigidity 
implies that the number of pairwise non-isomorphic $k$-graphs $F'$
with $n$ vertices that we can obtain this way is at least 
${n-v(F)+\ell-1\choose \ell-1}$ (the number of solutions to
$n-v(F)=x_1+\dots+x_{\ell}$ 
in non-negative integers) divided by $\ell!$. Moreover, each such
$F'$ will appear an an induced subgraph in every large
maximum $P$-construction by Lemmas~\ref{lm:regular} and~\ref{lm:OmegaN}.
Since $\ell$ can be
chosen arbitrarily large, $P$ is indeed complex.
Thus Theorem~\ref{th:Fin} answers 
the question
of Falgas-Ravry and Vaughan~\cite[Question~4.4]{falgas+vaughan:13}
to solve an explicit Tur\'an problem with a complex extremal
configuration (if one agrees that the family $\C F$ in Theorem~\ref{th:Fin}
is ``explicit'').

Let $\PI{k}{m}$ consist of all possible Tur\'an densities $\pi(\C F)$ 
where $\C F$ is a family consisting
of at most $m$ forbidden $k$-graphs. 

\begin{question}[Baber and Talbot~\cite{baber+talbot:12}]
 \label{qn:BT} Let $k\ge 3$. Which of the following trivial inclusions 
 $
 \PI{k}{1}\subseteq \PI{k}{2}\subseteq \dots \subseteq \PI{k}{i}
\subseteq\dots \subseteq \FinPi{k}
 $
 are strict?
\end{question}

It is open even whether $\PI{k}{1}=\FinPi{k}$ for $k\ge 3$.

\begin{question}[Jacob Fox (personal communication)]
Does $\FinPi{k}$ contain a transcendental number?\end{question}

Since there are only countably many algebraic numbers,  Theorem~\ref{th:uncountable} implies that $\InfPi{k}$ has a transcendental
number
for every $k\ge 3$.

\begin{question}[Frank, Peng, R\"odl, and Talbot~\cite{frankl+peng+rodl+talbot:07}]
Let $k\ge 3$. Is there $\alpha_k<1$ such that no value in  
$(\alpha_k,1)$ is a jump for $k$-graphs?
\end{question}

Note that by Proposition~\ref{pr:closed} the last condition is 
equivalent to $\InfPi{k}\supseteq [\alpha,1]$.
It is still open if $\InfPi{k}$ contains 
some interval of positive length for $k\ge 3$. 
On the other hand, the arsenal of 
tools for proving that some real
does \emph{not} belong to $\InfPi{k}$ is very limited for $k\ge 3$. 
In addition to the old result of 
Erd\H os~\cite{erdos:64}
that $\InfPi{k}\cap (0,k!/k^k)=\emptyset$, the only
other such result is by Baber and Talbot~\cite{baber+talbot:11} that 
$(0.2299,0.2315)\cap \InfPi{3}=\emptyset$. The proof
in~\cite{baber+talbot:11} uses flag algebras and is
computer-generated.

Hatami and Norine~\cite{hatami+norine:11} showed that the question
whether a given linear inequality in subgraph densities is always valid
is undecidable.

\begin{question}\label{qu:dec} Is the validity of $\pi(\C F)\le \alpha$
decidable, 
where the input is a finite family $\C F$ of $k$-graphs and a
rational number $\alpha$?
\end{question}

A related open question is whether every true inequality $\pi(\C F)\le \alpha$
admits a finite proof in Razborov's Cauchy-Schwarz calculus~\cite{razborov:07,razborov:10}
(see also~\cite[Appendix~A]{hatami+norine:11}).

If $k=2$, then the answer to Question~\ref{qu:dec} is
in the affirmative by the Erd\H os-Stone-Simonovits 
Theorem~\cite{erdos+stone:46,erdos+simonovits:66}. Brown,
Erd\H os, and Simonovits~\cite{brown+erdos+simonovits:85} 
obtained a positive solution to the version of Question~\ref{qu:dec} for
the class of directed multigraphs.

As we have already mentioned, our proof of Theorem~\ref{th:Fin} relies on 
the Strong Removal Lemma.
So the size of
the obtained family $\C F$ is huge (even for 
small concrete $P$). This is in contrast to many previous results
and conjectures that
forbid very few hypergraphs. The main place in our proof that makes
$|\C F|$ huge is the application of the Removal Lemma
in the proof of Lemma~\ref{lm:edit}. If, for some concrete $P$, Lemma~\ref{lm:edit} can be deduced
in an alternative way, then one might 
be able to obtain an explicit and reasonably sized $\C F$ for which Theorem~\ref{th:Fin} holds for all large $n$. (Note
that we did not try to optimise our other lemmas
for the sake of brevity and generality.)
So, some of our results and techniques might be useful for small forbidden families
as well. Also, the new ideas introduced for
proving Theorem~\ref{th:Fin} (in particular the
method of Lemma~\ref{lm:Max}) might be applicable to other 
instances of the Tur\'an problem.

\section*{Acknowledgements}

The author thanks Zolt\'an F\"uredi for helpful discussions and the anonymous
referee for the comments that greatly
improved the presentation of this paper.

\renewcommand{\baselinestretch}{0.98}\small


\begin{thebibliography}{10}

\bibitem{baber+talbot:11}
R.~Baber and J.~Talbot.
\newblock Hypergraphs do jump.
\newblock {\em Combin.\ Probab.\ Computing}, 20:161--171, 2011.

\bibitem{baber+talbot:12}
R.~Baber and J.~Talbot.
\newblock New {Tur\'an} densities for 3-graphs.
\newblock {\em Electronic J.\ Combin.}, 19:21pp., 2012.

\bibitem{balogh:02}
J.~Balogh.
\newblock The {Tur\'an} density of triple systems is not principal.
\newblock {\em J.\ Combin.\ Theory\ {\rm (A)}}, 100:176--180, 2002.

\bibitem{bollobas:74}
B.~{Bollob\'as}.
\newblock Three-graphs without two triples whose symmetric difference is
  contained in a third.
\newblock {\em Discrete Math.}, 8:21--24, 1974.

\bibitem{BCLSV:08}
C.~Borgs, J.~Chayes, L.~Lov{\'a}sz, V.~T. S{\'o}s, and K.~Vesztergombi.
\newblock Convergent sequences of dense graphs {I}: {Subgraph} frequencies,
  metric properties and testing.
\newblock {\em Adv.\ Math.}, 219:1801--1851, 2008.

\bibitem{brown+erdos+simonovits:85}
W.~G. Brown, P.~{Erd\H os}, and M.~Simonovits.
\newblock Algorithmic solution of extremal digraph problems.
\newblock {\em Trans.\ Amer.\ Math.\ Soc.}, 292:421--449, 1985.

\bibitem{brown+erdos+simonovits:84}
W.~G. Brown, P.~Erd{\H{o}}s, and M.~Simonovits.
\newblock Inverse extremal digraph problems.
\newblock In {\em Finite and infinite sets, {V}ol.\ {I}, {II} ({E}ger, 1981)},
  volume~37 of {\em Colloq. Math. Soc. J\'anos Bolyai}, pages 119--156.
  North-Holland, Amsterdam, 1984.

\bibitem{brown+simonovits:84}
W.~G. Brown and M.~Simonovits.
\newblock Digraph extremal problems, hypergraph extremal problems and the
  densities of graph structures.
\newblock {\em Discrete Math.}, 48:147--162, 1984.

\bibitem{brown+simonovits:02}
W.~G. Brown and M.~Simonovits.
\newblock Extremal multigraph and digraph problems.
\newblock In {\em Paul {E}rd{\H o}s and his mathematics, {II} ({B}udapest,
  1999)}, volume~11 of {\em Bolyai Soc. Math. Stud.}, pages 157--203. J\'anos
  Bolyai Math. Soc., Budapest, 2002.

\bibitem{chung+graham:eg}
F.~Chung and R.~L. Graham.
\newblock {\em {Erd\H{o}s} on Graphs: His Legacy of Unsolved Problems}.
\newblock A.K.Peters, Wellesley, 1998.

\bibitem{decaen+furedi:00}
D.~de~Caen and Z.~F{\"u}redi.
\newblock The maximum size of 3-uniform hypergraphs not containing a {F}ano
  plane.
\newblock {\em J.\ Combin.\ Theory\ {\rm (B)}}, 78:274--276, 2000.

\bibitem{decaen+kreher+wiseman:88}
D.~de~Caen, D.~L. Kreher, and J.~Wiseman.
\newblock On constructive upper bounds for the {T}ur\'an numbers
  {$T(n,2r+1,2r)$}.
\newblock {\em Congr. Numer.}, 65:277--280, 1988.

\bibitem{elek+szegedy:12}
G.~Elek and B.~Szegedy.
\newblock A measure-theoretic approach to the theory of dense hypergraphs.
\newblock {\em Adv.\ Math.}, 231:1731--1772, 2012.

\bibitem{erdos:64}
P.~Erd{\H{o}}s.
\newblock On extremal problems of graphs and generalized graphs.
\newblock {\em Israel J.\ Math.}, 2:183--190, 1964.

\bibitem{erdos:67a}
P.~Erd{\H{o}}s.
\newblock Some recent results on extremal problems in graph theory. {R}esults.
\newblock In {\em Theory of Graphs (Internat. Sympos., Rome, 1966)}, pages
  117--123 (English); pp. 124--130 (French). Gordon and Breach, New York, 1967.

\bibitem{erdos+simonovits:66}
P.~Erd{\H{o}}s and M.~Simonovits.
\newblock A limit theorem in graph theory.
\newblock {\em Stud.\ Sci.\ Math.\ Hungar.}, pages 51--57, 1966.

\bibitem{erdos+stone:46}
P.~Erd{\H o}s and A.~H. Stone.
\newblock On the structure of linear graphs.
\newblock {\em Bull.\ Amer.\ Math.\ Soc.}, 52:1087--1091, 1946.

\bibitem{falgas+vaughan:13}
V.~Falgas-Ravry and E.~R. Vaughan.
\newblock Applications of the semi-definite method to the {Tur\'an} density
  problem for {$3$}-graphs.
\newblock {\em Combin.\ Probab.\ Computing}, 22:21--54, 2013.

\bibitem{frankl+furedi:89}
P.~Frankl and Z.~F{\"u}redi.
\newblock Extremal problems whose solutions are the blowups of the small
  {W}itt-designs.
\newblock {\em J.\ Combin.\ Theory\ {\rm (A)}}, 52:129--147, 1989.

\bibitem{frankl+peng+rodl+talbot:07}
P.~Frankl, Y.~Peng, V.~R{\"o}dl, and J.~Talbot.
\newblock A note on the jumping constant conjecture of {E}rd{\H o}s.
\newblock {\em J.\ Combin.\ Theory\ {\rm (B)}}, 97:204--216, 2007.

\bibitem{frankl+rodl:84}
P.~Frankl and V.~R{\"o}dl.
\newblock Hypergraphs do not jump.
\newblock {\em Combinatorica}, 4:149--159, 1984.

\bibitem{furedi:91}
Z.~F{\"u}redi.
\newblock {Tur\'an} type problems.
\newblock In {\em Surveys in Combinatorics}, volume 166 of {\em London Math.\
  Soc.\ Lecture Notes Ser.}, pages 253--300. Cambridge Univ.\ Press, 1991.

\bibitem{furedi+pikhurko+simonovits:03:ejc}
Z.~F{\"u}redi, O.~Pikhurko, and M.~Simonovits.
\newblock The {Tur{\'a}n} density of the hypergraph {$\{abc,ade,bde,cde\}$}.
\newblock {\em Electronic J.\ Combin.}, 10:7pp., 2003.

\bibitem{giraud:90}
G.~R. Giraud.
\newblock Remarques sur deux probl\`emes extr\'emaux.
\newblock {\em Discrete Math.}, 84:319--321, 1990.

\bibitem{hatami+norine:11}
H.~Hatami and S.~Norine.
\newblock Undecidability of linear inequalities in graph homomorphism
  densities.
\newblock {\em J.\ Amer.\ Math.\ Soc.}, 24:547--565, 2011.

\bibitem{katona+nemetz+simonovits:64}
G.~O.~H. Katona, T.~Nemetz, and M.~Simonovits.
\newblock On a graph problem of {Tur\'an} {(In Hungarian)}.
\newblock {\em Mat.\ Fiz.\ Lapok}, 15:228--238, 1964.

\bibitem{keevash:11}
P.~Keevash.
\newblock Hypergraph {Tur\'an} problem.
\newblock In R.~Chapman, editor, {\em Surveys in Combinatorics}, pages 83--140.
  Cambridge Univ.\ Press, 2011.

\bibitem{lovasz:lngl}
L.~{Lov\'asz}.
\newblock {\em Large Networks and Graph Limits}.
\newblock Colloquium Publications. Amer.\ Math.\ Soc, 2012.

\bibitem{lovasz+szegedy:06}
L.~Lov{\'a}sz and B.~Szegedy.
\newblock Limits of dense graph sequences.
\newblock {\em J.\ Combin.\ Theory\ {\rm (B)}}, 96:933--957, 2006.

\bibitem{motzkin+straus:65}
T.~S. Motzkin and E.~G. Straus.
\newblock Maxima for graphs and a new proof of a theorem of {T}ur\'an.
\newblock {\em Can.\ J.\ Math.}, 17:533--540, 1965.

\bibitem{mubayi:06}
D.~Mubayi.
\newblock A hypergraph extension of {Tur\'an's} theorem.
\newblock {\em J.\ Combin.\ Theory\ {\rm (B)}}, 96:122--134, 2006.

\bibitem{mubayi+pikhurko:08}
D.~Mubayi and O.~Pikhurko.
\newblock Constructions of non-principal families in extremal hypergraph
  theory.
\newblock {\em Discrete Math.}, 308:4430--4434, 2008.

\bibitem{peng+zhao:08}
Y.~Peng and C.~Zhao.
\newblock Generating non-jumping numbers recursively.
\newblock {\em Discrete Applied Math.}, 156:1856--1864, 2008.

\bibitem{razborov:07}
A.~Razborov.
\newblock Flag algebras.
\newblock {\em J.\ Symb.\ Logic}, 72:1239--1282, 2007.

\bibitem{razborov:10}
A.~Razborov.
\newblock On 3-hypergraphs with forbidden 4-vertex configurations.
\newblock {\em {SIAM} J.\ Discr.\ Math.}, 24:946--963, 2010.

\bibitem{razborov:11:psim}
A.~A. Razborov.
\newblock On the {F}on-der-{F}laass interpretation of extremal examples for
  {Tur\'an's} {$(3,4)$}-problem.
\newblock {\em Proc.\ Steklov Inst.\ Math.}, 274:247--266, 2011.
\newblock Translated from \emph{Trudy Mat. Inst. Steklova}.

\bibitem{rodl+schacht:09}
V.~R{\"o}dl and M.~Schacht.
\newblock Generalizations of the {Removal} {Lemma}.
\newblock {\em Combinatorica}, 29:467--501, 2009.

\bibitem{ruzsa+szemeredi:78}
I.~Z. Ruzsa and E.~{Szemer\'edi}.
\newblock Triple systems with no six points carrying three triangles.
\newblock In A.~Hajnal and V.~{S\'os}, editors, {\em Combinatorics {II}}, pages
  939--945. North Holland, Amsterdam, 1978.

\bibitem{sidorenko:87}
A.~Sidorenko.
\newblock The maximal number of edges in a homogeneous hypergraph containing no
  prohibited subgraphs.
\newblock {\em Math Notes}, 41:247--259, 1987.
\newblock Translated from \emph{Mat.\ Zametki}.

\bibitem{sidorenko:93:c}
A.~Sidorenko.
\newblock Boundedness of optimal matrices in extremal multigraph and digraph
  problems.
\newblock {\em Combinatorica}, 13:109--120, 1993.

\bibitem{sidorenko:95}
A.~Sidorenko.
\newblock What we know and what we do not know about {Tur\'an} numbers.
\newblock {\em Graphs Combin.}, 11:179--199, 1995.

\bibitem{simonovits:68}
M.~Simonovits.
\newblock A method for solving extremal problems in graph theory, stability
  problems.
\newblock In {\em Theory of Graphs (Proc. Colloq., Tihany, 1966)}, pages
  279--319. Academic Press, 1968.

\bibitem{simonovits:74}
M.~Simonovits.
\newblock Extermal graph problems with symmetrical extremal graphs.
  {A}dditional chromatic conditions.
\newblock {\em Discrete Math.}, 7:349--376, 1974.

\bibitem{simonovits:83:misc}
M.~Simonovits.
\newblock Extremal graph problems and graph products.
\newblock In {\em Studies in pure mathematics}, pages 669--680. Birkh\"auser,
  Basel, 1983.

\bibitem{turan:41}
P.~{Tur\'an}.
\newblock On an extremal problem in graph theory (in {Hungarian}).
\newblock {\em Mat.\ Fiz.\ Lapok}, 48:436--452, 1941.

\end{thebibliography}

\end{document}